\documentclass[12pt]{article}

\usepackage[a4paper, total={6.4in, 8.5in}]{geometry}
\usepackage{amsmath,amssymb,latexsym,amsfonts}
\usepackage{stmaryrd}
\usepackage{shapepar,bm}
\usepackage[dvips,arrow,matrix,ps,color,line]{xy}
\usepackage{theorem,amsfonts,amssymb,amscd}
\usepackage{geometry}
\usepackage{makeidx}
\usepackage{eulervm}

\usepackage{graphicx}
\usepackage{graphics}

\usepackage{epstopdf}

\usepackage{amssymb,graphics,hyperref}

\hypersetup{
    colorlinks = true,
    linkcolor = blue,
    anchorcolor = red,
   citecolor = blue,
    filecolor = red,
    pagecolor = red,
    urlcolor = blue}

\newcommand{\mathsym}[1]{{}}

\usepackage[usenames]{color}
\definecolor{MyLightMagenta}{cmyk}{0.1,0.8,0,0.1}
\definecolor{MyDarkBlue}{rgb}{0.1,0,0.3}
\makeindex

\def\wb{[\bfb]}

\def\wbet{[\bm\beta]}

\def\bfz{{\mathbf z}}

\def\ovsig{\overline{\sigma}}

\def\bbeta{{\bm\beta}}
\def\bfb{{\mathbf b}}

\def\bfu{{\mathbf u}}

\def\bfw{{\mathbf w}}

\def\Ccal{{\mathcal C}}

\def\Vcal{\mathcal V}

\def\bfx{{\mathbf x}}

\def\sgn{\mathrm{sgn}}

\def\ZZ{\mathbb Z}

\def\Scal{{\mathcal S}}
\def\Ecal{{\mathcal E}}

\def\QQ{\mathbb Q}

\def\cocoa{{\hbox{\rm C\kern-.13em o\kern-.07em C\kern-.13em o\kern-.15em A}}}

\def\Dcal{\mathcal D}

\def\hom{{\mathrm{Hom}}}

\def\bfu{{\bf u}}

\def\ovDcal{\overline{\Dcal}}

\def\End{\mathrm{End}}

\def\blamb{{\bm \lambda}}
\def\bmu{{\bm \mu}}
\def\bnu{{\bm \nu}}
\def\Pcal{{\mathcal P}}

\def\w2M{\bigwedge^2M}

\def\w{\wedge }
\def\bw{\bigwedge }
\def\bwV{{\bigwedge\hskip-3.5pt V}}
\def\bwVz{{\bigwedge\hskip-3.5pt V}\llb z\rrb}
\protect
\protect
\protect
\protect

\protect
\protect

\def\sra{\rightarrow}

\def\proof{\noindent{\bf Proof.}\,\,}
\def\qed{{\hfill\vrule height4pt width4pt depth0pt}\medskip}
\def\be{\begin{equation}}
\def\ee{\end{equation}}
\def\bclm{\begin{claim}}
\def\eclm{\end{claim}}
\def\beqn{\begin{eqnarray}}
\def\eeqn{\end{eqnarray}}
\def\beqn*{\begin{eqnarray*}}
\def\eeqn*{\end{eqnarray*}}

\theoremstyle{change}
\theorembodyfont{\rmfamily}

\newtheorem{claim}{}[section]
\global
\global

\def\no@breaks#1{{\def\\{ \ignorespaces}#1}}    % disallow explicit line breaks

  \makeatletter
\def\cleardoublepage{\clearpage\if@twoside \ifodd\c@page\else
\hbox{} \thispagestyle{empty}
\newpage
\if@twocolumn\hbox{}\newpage\fi\fi\fi} \makeatother

\usepackage{eso-pic,graphicx}
\makeatletter
\newcommand\BackgroundPicture[2]{%
  \setlength{\unitlength}{1pt}%
  default \put(0,\strip@pt\paperheight){%
  \parbox[t][\paperheight]{\paperwidth}{%
    \vfill
     \centering \includegraphics[angle=#2, width=15cm, height=15cm,  bb=0 0 150 150]{#1}
    \vfill
}}} %
\makeatother
\newcommand{\llb}{\llbracket}
\newcommand{\rrb}{\rrbracket}
\def\bmtx{\begin{matrix}}
\def\emtx{\end{matrix}}
\makeindex

\begin{document}

\date{}
\title{Polynomial ring representations \\ of endomorphisms of exterior powers}

\author{Ommolbanin Behzad, Andr\'e Contiero, \\ 
Letterio Gatto \& Renato Vidal Martins
\thanks{Work sponsored by  Finanziamento 
Diffuso della Ricerca, no. 53$\_$RBA17GATLET del Politecnico di Torino; 
\, Progetto di Eccellenza\, Dipartimento\, di\, Scienze\, Matematiche, 2018--2022 no.
E11G18000350001, INDAM-GNSAGA e PRIN "Geometria delle Variet\`a Algebriche''.
\smallskip 
\newline  ${}$
\,\,\,\,\,\,\,{\em Keywords and Phrases:} Hasse-Schmidt Derivations and Vertex Operators  on Exterior Algebras,  Bosonic and Fermionic Representations by Date-Jimbo-Kashiwara-Miwa, Symmetric Functions.
\smallskip
\newline ${}$ \,\,\,\,\,\, {\bf 2020 MSC:} 14M15, 15A75, 05E05, 17B69. \,}}
%\\${}$\\ { %(Bulletin Math. Braz. Soc. New Series, 2020}\\ \href{https://doi.org/10.1007/s00574-020-00195-9}{\tt https://doi.org/10.1007/s00574-020-00195-9})}

\maketitle

\begin{abstract}
\noindent
%The fermionic Fock representation of the canonical Clifford algebra supported on the direct sum of a countable infinite dimensional vector space  with its restricted dual implies  that a
A polynomial ring with rational coefficients is  an irreducible representation of  Lie 
algebras of  endomorphisms of  exterior powers of a  infinite countable dimensional $\QQ$--vector space. We give an explicit description of it, using suitable vertex operators on exterior algebras, which mimick those occurring in the bosonic vertex representation of the Lie algebra $gl_\infty$,  due to Date--Jimbo--Kashiwara and Miwa (DJKM).
\end{abstract}

\tableofcontents

\section{Introduction}

\claim{\bf The goal.} Let $B_r:=\QQ[e_1,\ldots,e_r]$ be the polynomial ring in $r$ indeterminates $(e_1,\ldots,e_r)$. This paper supplies its  explicit description  as a module over the
 Lie algebras of endomorphisms of $k$-th exterior powers of a vector 
space $V$ of infinite countable dimension.  The goal is achieved by means of  certain 
{\em vertex operators on the exterior algebra}, defined by means of {\em Schubert 
derivations}. 

The latter are distinguished Hasse-Schmidt derivations on exterior 
algebras, introduced in \cite{SCHSD} and extensively treated in \cite{HSDGA}; see also the survey \cite{Nigran} or  \cite[p.~116]{CoDuGu}, for more related discussions. They have 
shown their versatility in applications to improve effectiveness in Schubert 
Calculus computations (see \cite{NBFSC, tango}), to equivariant cohomology of 
Grassmannians (Cf. \cite{ESC}, but also \cite{LakEq1,LakEq2}), to generalise the Cayley-Hamilton theorem \cite{GaRow,GSCH} or, 
like in \cite{gln,SDIWP} and in the present paper, to revisit  the bosonic vertex 
representation of Lie algebras of endomorphisms as in \cite{DJKM01} (see also 
\cite{jimbomiwa} and \cite[Propositions 5.2--5.3]{KacRaRoz}), providing new methods and new insight. 

The 
$gl(\bw^kV)$--module structure of $B_r$, that we are going to describe, will be referred to as {\em bosonic representation} of 
$\bw^kV$, by a possibly strong, but suggestive, abuse of terminology, due to the evident relationship with pioneering work by Date, Jimbo, Kashiwara and 
Miwa \cite{DJKM01} (see 
also \cite{jimbomiwa} and \cite{KacRaRoz}) and a more general framework that, in \cite{Cox-Futorny}, one refers to as DJKM (affine) Heisenberg algebra.

That $B_r$ is a representation of $gl(\bw^kV)$  is easy to see in very special cases. For $k=0$, it is just 
multiplication by rational numbers, as $\bw^0V=\QQ$. For $k>r$,  is the trivial null 
representation. For $r=k=1$,  it amounts to the well known
general fact that any vector space is a module over the Lie algebra of its own endomorphisms. Thus,  the linear extension of the set map $e_1^i\mapsto b_i$ is a vector space 
isomorphism $B_1\sra V$, making $B_1$ into a $gl(V)$-module, by  pulling back 
that structure from $V$. For the general case see below, in the second  part of this introduction, devoted to state precisely the main result. The third one will discuss, instead, background and motivation.

\claim{\bf Statement of the main result.} The  ring $B_r$ possesses a $\QQ$-basis formed by certain Schur 
determinants $\Delta_\blamb(H_r)$ (like in Section~\ref{sec2:22}; see \cite[]{MacDonald}), where $\blamb$ ranges over the set $\Pcal_r$ of all the partitions of length 
at most $r$.  Let $V:=\bigoplus_{i\geq 0}\QQ\cdot b_i$ be the $\QQ$-vector space with 
basis $\bfb:=(b_i)_{i\geq 0}$. Then $\bw^rV:=\bigoplus_{\blamb\in\Pcal_r} \QQ\cdot \wb^r_\blamb$, where $
\wb^r_\blamb:=b_{r-1+\lambda_1}\w\cdots\w b_{\lambda_r}$. It follows that the linear 
extension of the set map $\Delta_\blamb(H_r)\mapsto\wb^r_\blamb$ is a $\QQ$-vector 
space isomorphism $B_r\sra \bw^rV$ sending $1\mapsto \wb^r_0:=b_{r-1}\w\cdots\w b_0$. It can be phrased by saying that $\bw^rV$ carries a structure of free $B_r$-module of rank $1$ generated by $
\wb^r_0$, such that $\wb^r_\blamb=\Delta_\blamb(H_r)\wb^r_0$ (Section \ref{sec2:22}).
%The other, 
%more refined, is based on the result of Laksov and Thorup who show that the $r$-the 
%exterior power $\bw^r\QQ[X]$ of the polynomial ring $\QQ[X]$ is a free module of rank 
%$1$ over the ring of symmetric functions $\QQ[X_1,\ldots, X_r]^{S_r}$ which, by the 
%fundamental theorem of symmetric functions, is itself a polynomial ring in the 
%elementary symmetric polynomials in $(X_1,\ldots, X_r)$, which is what motivate our 
%notation $B_r=\QQ[e_1,\ldots, e_r]$.

The restricted dual of $V$ is $V^*:=\bigoplus_{j\geq 0} \QQ\cdot\beta_j$, where  $\beta_j:V\sra \QQ$ is the unique linear form $\beta_j(b_i)=\delta_{ji}$. Then $\bw^rV^*:=\bigoplus_{\blamb
\in\Pcal_r}\QQ\cdot\wbet^r_\blamb$, where   $\wbet^r_\blamb(\wb^r_\bmu)=\delta_{\blamb,\bmu}$. Let $gl(\bw^kV)$ be the Lie algebra
of the endomorphisms of $\bw^kV$ vanishing at  $
\wb^k_\blamb$ for all partitions $\blamb\in \Pcal_k$ but finitely many. If  $\Ecal^k_{\bmu,\bnu}:=
\wb^k_\bmu\otimes\wbet^k_\bnu$, then: 
$$
gl(\bw^kV)=\bw^kV\otimes \bw^kV^*=
\bigoplus_{\bmu,\bnu\in\Pcal_k}\QQ\cdot \Ecal^k_{\bmu,\bnu}.$$ 

\noindent
The  $B_r$ (bosonic) representation of  $gl(\bw^kV)$,  for all $k,r\geq 0$, is then naturally defined via the following equality:
\be
(\Ecal^k_{\bmu,\bnu}\Delta_\blamb(H_r))\wb^r_0=\wb^k_\bmu\w (\wbet^k_\bnu\lrcorner
\wb^r_\blamb),\label{eq0:defrep}
\ee
where the contraction $\wbet^k_\bnu\lrcorner$ maps $\bw^rV$ to $\bw^{r-k}V$ (Section \ref{pairing}). To express the  $gl(\bw^kV)$--action \eqref{eq0:defrep} on $B_r$ within a compact formula, a standard philosophy suggests to use generating functions. Let $\bfz_k:=(z_1,\ldots, z_k)$ and $\bfw_k:=(w_1,
\ldots,w_k)$ be two sets of formal variables. The $k$-tuples of the formal inverses $
(z_1^{-1},\ldots,z_k^{-1})$ and $(w_1^{-1},\ldots, w_k^{-1})$   will be 
denoted by $\bfz_k^{-1}$  and $\bfw_k^{-1}$ respectively.  The standard notation
$s_\bmu(\bfz_k)$ and $s_\bnu(\bfw_k^{-1})$ stands for the symmetric Schur polynomials in the 
variables $\bfz_k$ and $\bfw_k^{-1}$ (See \cite[p.~40]{MacDonald}). Define 
$$
\Ecal(\bfz_k,\bfw_k^{-1})=\sum_{\bmu,\bnu\in\Pcal_k}\Ecal^k_{\bmu\bnu}\cdot s_
\bmu(\bfz_k)s_\bnu(\bfw_k^{-1}):B_r\sra B_r\llb \bfz_k,\bfw_k^{-1}]
$$ 
through
$$
\Big(\Ecal(\bfz_k,\bfw_k^{-1})\Delta_\blamb(H_r)\Big)\wb^r_0=\sum_{\bmu.
\bnu\in\Pcal_k}s_\bmu(\bfz_k)s_\bnu(\bfw_k^{-1})\wb^k_\blamb\w \left(\wbet^k_\bnu
\lrcorner \wb^r_\blamb\right).
$$
Our main result is:
\paragraph{\bf Theorem \ref{thm:mainthm}.} {\em The equality below holds for all $k,r\geq 0$ and all $\blamb\in\Pcal_r$:
\be
(\Ecal(\bfz_k,\bfw_k^{-1})\Delta_\blamb(H_r))\wb^r_0=\prod_{j=1}^k\left({z_j\over w_j}
\right)^{\hskip-3pt r-k}\cdot \Gamma(\bfz_k)\Gamma^*(\bfw_k)\wb^r_\blamb,\label{eq0:forthm}
\ee
}

\noindent
where the vertex operators  $\Gamma(\bfz_k),\Gamma^*(\bfz_k):\bwV\sra \bwV\llb\bfz_k,\bfz_k^{-1}]$  on the exterior algebra $\bwV$ of $V$ are  introduced in 
Definition~\ref{def4:vertex} and studied in more details in Sections \ref{sec:sec5} and 
\ref{sec:sec6}. They are merely defined  as product of {\em Schubert derivations}.  In case $r$ is big with respect to the length of $\blamb$, the vertex operators involved in formula \eqref{eq0:forthm} can be expressed as
$$
\Gamma(\bfz_k):=\prod_{j=1}^k{1\over E_r(z_j)}\exp\left(-\sum_{i\geq 1}{1\over i}\delta(\sigma_{-1}^i)p_i(\bfz_k^{-1})\right)
$$
and
$$
\Gamma^*(\bfw_k):=\prod_{j=1}^k{E_r(w_j)}\exp\left(\sum_{i\geq 1}{1\over i}\delta(\sigma_{-1}^i)p_i(\bfw_k^{-1})\right),
$$
where $E_r(z)$ is the generic monic polynomial $1-e_1z+\cdots+(-1)^re_rz^r$, the map  $\sigma_{-1}$ is the locally nilpotent endomorphism of $V$ mapping $b_j\mapsto b_{j-1}$ if $j\geq 1$ and $b_0\mapsto 0$, $\delta:gl(V)\mapsto \End(\bw V)$ is the natural representation of $gl(V)$ as a Lie algebra of (even) derivations of $\bw V$ and, finally,  $p_i(\bfu_k)$  denotes the Newton power sum $u_1^i+\cdots+u_k^i$ of degree $i$. 

In other words,  the image of $\Delta_\blamb(H_r)$ through $\Ecal^k_{\bmu,\bnu}$ is the coefficient of $s_\bmu(\bfz_k)s_\bnu(\bfw_k^{-1})$ in the right hand side of \eqref{eq0:forthm}. This may sounds tricky to evaluate, but it coincides with  the coefficient of 
$$
z_1^{k-1+\mu_1}\cdots z_k^{\mu_k}\cdot w_1^{-k+1-\nu_1}\cdots w_k^{-\nu_k}
$$
of the second member of \eqref{eq0:forthm}, multiplied by the Vandermonde determinants of $\bfz_k$ and $\bfw_k^{-1}$.

\claim{\bf Background and Motivations.} 
%Let $\Vcal:=\QQ[X^{-1},X]$ be the vector space of Laurent polynomials. One motivation for the present investigation  is that for $k=1$  formula  \eqref{eq0:forthm} is the shadow  of the limit case attained putting $r=\infty$. In this case  $B:=B_\infty$ would be 
%the polynomial ring in infinitely many indeterminates and our formula  is related with the 
%%celebrated DJKM bosonic vertex representation of $gl(\Vcal)$. That even this classical case can be treated  within the formalism of Schubert derivations has been shown  in the contribution \cite{SDIWP}.
%Each $\QQ$--vector space is endowed with an obvious module structure over the Lie sub-algebra  $V\otimes V^*\subseteq V$ of $End_\QQ(V)$.  The exterior algebra $\bw V$ of $V$ does make any excteption. It is then a representation of the Lie super-algebra algebra $gl(\bw V):=\bwV\otimes \bwV^*$. It is well known that $gl(\bw V)$ coincides with the canonical Clifford algebra $\Ccal:=\Ccal(V\oplus V^*)$ supported on $V\oplus V^*$, whose elements are finite linear combinations of finite words $b_{i_1}\cdots b_{i_h}\beta_{j_1}\cdots\beta_{j-k}$. Since the fermionic Fock space as in \cite{Frenkel-zvi} is itself a  principal module over the Clifford algebra $\Ccal$ freely generated by the vacuum vector $|0\rangle$ which is sometimes suggestively written as $b_0\w b_{-1}\w b_{-2}\w\cdots$.
 This paper is the first step towards the authors' attempt to better understand a fundamental, although elementary, representation theoretical fact. Let $
\Vcal:=\bigoplus_{j\in\ZZ}\QQ\cdot b_j$ be a vector space with basis $(b_j)$ parameterized 
by the integers (one may think of $\Vcal$ as being the vector space of the Laurent 
polynomials) and $\Vcal^*$ its restricted dual with basis $(\beta_j)_{j\in\ZZ}$. 

It is well known
that $\Vcal\oplus \Vcal^*$ supports a canonical structure of Clifford algebra $\Ccal:=
\Ccal(\Vcal\oplus\Vcal^*)$ \cite[p.~85]{Frenkel-zvi} or \cite{HSDGA} and that the {\em 
Fermionic Fock space} $F$ (also called the semi infinite wedge power $\bw^{\infty/2}\Vcal$, 
see \cite{BloOko01})  is an irreducible representation of $\Ccal$. 
More 
precisely, $F$ is an invertible module over the Lie super-algebra $\Ccal$ generated 
by a distinguished vector  $|0\rangle$, the {\em vacuum}, that in the 
formalism of the infinite wedge powers can be suggestively written as $b_0\w b_{-1}\w 
b_{-2}\w\cdots$ The huge Clifford algebra $\Ccal$, whose elements are finite 
linear combinations of words of the form $b_{i_i}\cdots b_{i_h}\beta_{j_1}\cdots
\beta_{j_k}$, contains in a natural way all, but not only,  the  Lie algebras 
$gl(\bw^kV)$, for all $k\geq 0$. In particular, it turns out that $F$ is a $gl(\bw^k
\Vcal)$-module for all $k\geq 0$. Then, the bosonic Fock space  $B:=B_\infty:=
\QQ[e_1,e_2,\ldots]$ gets a $gl(\bw^kV)$-module structure, for all $k\geq 0$, pulling 
back that of $F$ via the {\em boson-fermion correspondence}, a natural module 
isomorphism $B\sra F$  over the infinite dimensional Lie Heisenberg algebra. 

The latter may well be interpreted as a sort of Poincar\'e duality for infinite dimensional 
 Grassmannians. This case shall be analyzed in a forthcoming paper: although the formal 
 framework looks the same, the case $r=\infty$ is not just a naive limit of our formula 
 \eqref{eq0:forthm}, as indicated by the presence of the factor $\prod(z_j/w_j)^{r-k}$. 
 Since the algebra of endomorphisms of the exterior algebra of $V\cong\QQ[X]$ is 
 precisely the same Clifford algebra $\Ccal$ we alluded to above, we realised that was 
 already relevant and interesting to give a first closer look  to the $gl(\bwV)$--structure of $\bwV$, certainly not treated in any literature we have consulted up to now. The task, however, does not 
 look easy, at first sight, because $gl(\bwV)$ also contains the vector spaces  $\hom_
 \QQ(\bw^{k_1}V, \bw^{k_2}V)\cong \bw^{k_2}V\otimes \bw^{k_1}V^*$, with $k_1\neq k_2$. Thus, in this paper we are going to offer the description of the easiest case, namely the  representation 
 of homogeneous endomorphisms of $\bwV$ of degree $0$ (with respect to the exterior algebra 
 graduation). 
 
 The output is that the direct sum $\bigoplus_{k\geq 0}gl(\bw^kV)$ is 
 a Lie subalgebra of $gl(\bwV)$, represented by $B_r$ for all $r\geq 0$.
In the case of the fermionic Fock space, the $gl(\bw^{\hskip-2pt1}{\hskip-2pt V})$-structure of $B_\infty$ is the 
DJKM one \cite{DJKM01,jimbomiwa}. The reference \cite{SDIWP} already shows how the techniques of this paper 
are also suited to cope with the DJKM case and we know that the same methods will work 
as well for the DJKM representation of $gl(\bw^k\Vcal)$. To our knowledge no such a 
computation is known for the case $k>1$, and to gain feeling and experience we first coped with the more classical problem of describing the $gl(\bwV)$-module 
structure of $\bwV$. As a byproduct, examples of many computations with vertex operators occurring in the Heisenberg vertex algebra have been provided. Indeed, the references  \cite{HSDGA,pluckercone,SDIWP} already show our vertex operators on Grassmann algebras tend, as $r\sra \infty$ to the well known classical ones as in, e.g., \cite[Theorem 5.1]{KacRaRoz}.
\claim{\bf Organisation of the paper.} 
%The paper, intended to be as a first draft for submission, will probably exceed in %providing pre-requisites in the attempt to be as much self contained as possible. 
Most of preliminaries and basic notation are exposed in Section \ref{sec1:notation}. 
To be as much self contained as possible, the first part recalls  basics borrowed from the elementary theory of symmetric 
polynomials, such as, e.g. in \cite{MacDonald}, while  a second part quickly accounts  
on the notion of  Hasse-Schmidt derivation on Exterior algebras. The notion of inverse and transpose, essential for the sequel, are also discussed. 

The very special case of  HS derivations we will be concerned with are those that already in the very first reference \cite{SCHSD} was termed Schubert derivation. 
In a finite dimensional context,  the Schubert derivation is 
nothing but the Chern polynomial of the universal quotient bundle over a Grassmannian: it 
acts as a HS--derivation on the exterior algebra of the homology of the projective 
space, which is the same  as saying that to dealing with Schubert calculus for Grassmannians, B\'ezout theorem suffices. 

The explicit expression of the Schubert derivations offered in 
Section~\ref{sec:sec3}, introduced for the first time in this paper, makes all evident their strict relationship 
with vertex operators.  The Schubert derivations we consider are denoted by $\sigma_+(z),\ovsig_+(z)$, $\sigma_-(w)$ and $\ovsig_-(w)$. Those with the same sign as subscripts commute in the algebra of endomorphisms of the exterior algebra. However, due to the fact that $\sigma_-(w)$ and $\ovsig_-(w)$ are locally nilpotent, they commute with $\sigma_+(z)$ and $\ovsig_+(z)$ just up to the multiplication by a rational function.
%We consider four Schubert derivations, two of them being locally nilpotent, which causes that their mutual product,  a (multivariate) HS-derivation is not  but they do commute only up the 
%multiplication to a rational function which is related with the structure of a 
%Vandermonde determinant. 
To fully check our main Theorem \ref{thm:mainthm}, we devote Section \ref{sec:sec4}  to state and 
prove some relevant  commutation rules, some of which can be recognized within the phrasing of the  categorical framework for the Boson-Fermion correspondence, depicted in  \cite{FrPeSe}. See also \cite{Tian01} for a recent update.  

Section~\ref{sec:sec4} also contains an elegant definition of what we have proposed to name {\em vertex operators on a Grassmann 
algebra}, because its  evident relationship with those occurring in the classical boson--fermion 
correspondence. They arise, however,  in a definitely more elementary context that, in our opinion, 
would deserve to be further investigated, as we are currently doing. 

Vertex operators in 
the sense of Definition \ref{def4:vertex} are homogeneous operators of the exterior 
algebra, one of positive and the other of negative degree. We devote one section to 
each one of them (Sections~\ref{sec:sec5} and \ref{sec:sec6}) to dig up their 
relationship with basic computations in multi-linear algebra, such as wedging and 
contracting. Of course this idea is already present in the infinite wedge power 
context (e.g. \cite[Chapter 5]{KacRaRoz}), but this paper, together with \cite{pluckercone,gln,SDIWP}, is the first 
instance of applications of those techniques and ideas in finite dimensional landscapes. 
Finally, last Section~\ref{sec:sec7} is concerned with the statement and the proof of the main theorem togethr with some of its straightforward declinations in terms of familiar objects, like suitable deformations of the same Giambelli's determinants occurring in classical Schubert Calculus -- see Theorem~\ref{thm:them79}. To achieve the proof of the main Theorem~\ref{thm:mainthm}, we needed some preliminary Lemmas (such as \ref{lem:genbl} and \ref{lem:genbet}) which, along with Theorems \ref{thm:thm61} and \ref{thm:thm71}, we believe interesting in their own, as pieces of multilinear algebra properties addressed to wider general mathematical audiences.

\section{Preliminaries and Notation}\label{sec1:notation}

\claim{} The main characters of this paper are polynomials rings $B_r:=\QQ[e_1,
\ldots,e_r]$ in the $r\geq 1$ indeterminates $(e_1,\ldots,e_r)$ and a $\QQ$-vector 
space $V:=\bigoplus_{i\geq 0}\QQ\cdot b_i$ with basis
$\bfb:=(b_i)_{i\geq 0}$. The {\em restricted dual of} $V$ is $V^*:=\bigoplus_{j\geq 0}\QQ\cdot
\beta_j$,   where $\beta_j(b_i)=\delta_{ij}$.
We denote by $\bfb(z)$ and $\bbeta(w^{-1})$ the generating series of the  basis 
elements of $V$ and of $V^*$, i.e.:
\be
\bfb(z):=\sum_{i\geq 0}b_iz^i\qquad \mathrm{and}\qquad \bbeta(w^{-1}):=\sum_{j\geq 0}
\beta_jw^{-j}.\label{eq2:generating}
\ee
The exterior algebra of $V$ is $\bwV:=\bigoplus_{j\geq 0}\bw^jV$, the direct sum of the 
exterior powers $\bw^jV$, where $\bw^0V=\QQ$ and $\bw^1V=V$. The algebra structure is 
given by the $\QQ$-linear extension of the juxtaposition.
Given the generic polynomial $E_r(z):=1-e_1z+\cdots+(-1)^r e_r z^r\in B_r[z]$, one 
defines the sequence $H_r:=(h_j)_{j\in\ZZ}$ through the equality:
\be
\sum_{n\in\ZZ}h_nz^n:={1\over E_r(z)}.\label{eq:defhj}
\ee
In particular $h_j=0$ if $j<0$ and $h_0=1$. Moreover for $j\geq 0$, $h_j$ is an 
explicit polynomial homogeneous of degree $j$ in $(e_1,\ldots, e_r)$, once one gives 
weight $i$ to $e_i$.
  
\claim{} \label{sec2:22} A \emph{partition} is a monotonic non increasing sequence $
\blamb:=(\lambda_1\geq\lambda_2\geq\dots)$ 
of non negative integers, said to be its {\em parts}. Its length  $\ell(\blamb)$ is the 
number of its non zero parts, and $|\blamb|=\sum_{i\geq 0} \lambda_i$ is its {\em 
weight}. Let $\Pcal_r$ the set of all partitions of length at most $r$. To each $\blamb
\in\Pcal_r$ we associate 
\be
\wb^r_\blamb:=b_{r-1+\lambda_1}\wedge b_{r-2+\lambda_2}\wedge\dots\wedge b_{\lambda_r}
\in \bw^rV,
\ee
so that $(\wb^r_\blamb)_{\blamb\in\Pcal_r}$ is a $\QQ$-basis of $\bw^rV$,
and the {\em Schur determinant}
\be
\Delta_\blamb(H_r):=\det(h_{\lambda_j-j+i})_{1\leq i,j\leq r}=\begin{vmatrix}
h_{\lambda_1} & h_{\lambda_2-1} &\ldots & h_{\lambda_r -r+1} \\
h_{\lambda_1+1} & h_{\lambda_2} &\ldots & h_{\lambda_r-r+2} \\
%\vdots  & \vdots & \ddots & \vdots \\
%w_1^{-\lambda_1-k+1} & w_1^{-\lambda_2-k+2} &\ldots & w_1^{-\lambda_r+k-r} \\
%h_{\lambda_1+k} & h_{\lambda_2+k-1} &\ldots &  h_{\lambda_r+k-r+1} \\
\vdots  & \vdots & \ddots & \vdots \\
h_{\lambda_1+r-1} & h_{\lambda_2+r-2} &\ldots &  h_{\lambda_r} \\
\end{vmatrix}.\label{eq1:schurdet}
\ee
It is well known that
\be
B_r:=\bigoplus_{\blamb\in\Pcal_r}\QQ\cdot \Delta_\blamb(H_r)\label{eq1:modBr}
\ee
i.e. the Schur determinants form a $\QQ$-vector space basis of $B_r$ parametrized by 
the partitions of length at most $r$. It follows that $B_r$ is naturally isomorphic to 
$\bw^rV$ via the $\QQ$-linear extension of the sets map 
\be
\Delta_\blamb(H_r)\mapsto \wb^r_\blamb.\label{eq1:fbfc}
\ee

\claim{} Especially in the last section we shall be concerned with Schur polynomials in a set of indeterminates. We recall them here.
For each partition of length at most $k$ and any set of $k$ formal variables $\bfx_k:=(x_1,\ldots,x_k)$, one defines
$$
\Delta_\blamb(\bfx_k)=\det(x_j^{\lambda_{k-i+1}+i-1})=\begin{vmatrix} x_1^{\lambda_k}&x_2^{\lambda_k}&\cdots&x_k^{\lambda_k}\cr\cr
x_1^{1+\lambda_{k-1}}&x_2^{1+\lambda_{k-1}}&\cdots&x_k^{1+\lambda_{k-1}}\cr
\vdots&\vdots&\ddots&\vdots\cr
x_1^{k-1+\lambda_1}&x_2^{k-1+\lambda_1}&\cdots&x_k^{k-1+\lambda_1}
\end{vmatrix}.
$$
This is an skew symmetric polynomials in $(x_1,\ldots,x_k)$ and therefore divisible by the {\em Vandermonde determinant}
$$
\Delta_0(\bfx_k):=\begin{vmatrix} 1&1&\cdots& 1\cr\cr
x_1&x_2&\cdots&x_k\cr
\vdots&\vdots&\ddots&\vdots\cr 
x_1^{k-1}&x_2^{k-1}&\cdots&x_k^{k-1}
\end{vmatrix}=\prod_{i<j}(x_j-x_i).
$$
The Schur polynomial associated to $\bfx_k$ and the partition $\blamb$ is defined by the equality
$$
\Delta_\blamb(\bfx_k)=s_\blamb(\bfx_k)\cdot \Delta_0(\bfx_k),
$$
often said to be the {\em Jacobi-Trudy} formula.

\claim{}\label{pairing} The pairing
\be
(\beta_{i_1}\w\cdots\w\beta_{i_r})(v_1\w\cdots\w v_r)=\begin{vmatrix}\beta_{i_1}(v_1)&
\cdots&\beta_{i_1}(v_r)\cr\vdots&\ddots&\vdots\cr
\beta_{i_r}(v_1)&\cdots&\beta_{i_r}(v_r)
\end{vmatrix}\label{eq1:pairing}
\ee
establishes a natural identification between $\bw^rV^*$ and $(\bw^rV)^*$. If one 
denotes by $\wbet^r_\bmu$ the basis element 
\be 
\beta_{r-1+\mu_1}\w\cdots\w \beta_{\mu_r}\label{eq1:wbet}
\ee 
of $\bw^rV^*$, an easy check shows that  $\wbet^r_\bmu(\wb^r_\blamb)=\delta_{\bmu,
\blamb}$. The pairing \eqref{eq1:pairing} enables to attach to any $\beta\in V^*$ a map 
$\beta\lrcorner: \bwV\sra \bwV$ of degree $-1$ (with respect to the graduation of the 
exterior algebra) via the equality 
\be
\eta(\beta\lrcorner u)=(\beta\w\eta)(u),\qquad\quad \forall (u,\eta)\in \bw^rV\times 
\bw^{r-1}V^*.
\ee

\claim{} Let  $\bwV\llb z\rrb$ denote the formal power series in the indeterminate $z$ 
with  coefficients in the exterior algebra $\bw V$ of $V$. If  $\Scal$ is any 
set of indeterminates over $\QQ$, denote by $\QQ[\Scal]$ the corresponding algebra of 
formal power series.
The following is an extended reformulation of the main definition of the reference 
\cite{SCHSD} (see also \cite{HSDGA}).
\claim{\bf Definition.}\label{def:HSDr} 
 {\em  By {\em Hasse--Schmidt} derivation on $\bwV$ we mean any $\QQ[\Scal]$-linear 
 extension of a $\QQ$-linear map $\Dcal(z):\bwV\sra \bwV\llb z\rrb$
such that
\be
\Dcal(z)(u\w v)=\Dcal(z)u\w\Dcal(z)v,\qquad \forall u,v\in\bwV,\label{eq:eq2}
\ee
which, by abuse of notation, will be  denoted by the same symbol  
$$
\Dcal(z):\QQ\llb\Scal\rrb\otimes_\QQ\bwV\sra \QQ[\Scal]\otimes_\QQ\bwV\llb z\rrb,
$$
{\em (}instead of the more precise, but lengthier, $1_{\QQ\llb\Scal\rrb}\otimes_\QQ
\Dcal(z)${\em )}.}

\claim{} If $\Dcal_i\in\End_\QQ(\bwV)$ are such that
$
\sum_{i\geq 0}\Dcal_iz^i:=\Dcal(z).
$
then \eqref{eq:eq2} is equivalent to the  system of relations holding for all $i\geq 0$
$$
\Dcal_i(u\wedge v)=\sum_{j=0}^{i}\Dcal_{j}u\wedge D_{i-j}v,
$$ 
By \cite[Proposition 4.1.7, (iii)]{HSDGA}, if $\Dcal_0$ is 
invertible in $\End_\QQ(\bwV)$, then $\Dcal(z)$ is invertible as a $\End_\QQ(\bw V)$-
valued formal power series and its  inverse, $\ovDcal(z)$, is an HS--derivation as 
well. The {\em integration by parts formulas} follow for all $u,v\in\bwV$:
\begin{eqnarray}
\Dcal(z)(\ovDcal(z)u\w v)&=&u\w \Dcal(z)v,\label{eq:intp1}\\ \nonumber \\
\ovDcal(z)(\Dcal(z)u\w v)&=&u\w \ovDcal(z)v.\label{eq:intp2}
\end{eqnarray}
Formulas \eqref{eq:intp1} and \eqref{eq:intp2}  are implicitly assuming the $\QQ[[z]]$-
linearity of $\Dcal(z)$ we alluded to in Definition \ref{def:HSDr}. The extension of 
the linearity an HS-derivation over polynomials algebra will be assumed in the 
following without any further mention.

\claim{} The notation
\be
\Dcal(z)\wb^r_\blamb=\big[\Dcal(z)\bfb\big]^r_\blamb
\ee
will be used as a shorthand for the equality
$$
\Dcal(z)\wb^r_\blamb=\Dcal(z)\big(b_{r-1+\lambda_1}\w\cdots\w b_{\lambda_r}\big)=\Dcal(z)b_{r-1+\lambda_1}\w\cdots\w \Dcal(z) b_{\lambda_r}
$$
meaning that $\Dcal(z)$ is a HS-derivation.
\claim{} \label{sec:transp}The transpose  $\Dcal(z)^T:\bwV^*\sra\bwV^*\llb z\rrb$ of the HS derivation $
\Dcal(z)$ is defined via its action on homogeneous elements. If $\eta\in \bw^rV^*$, 
then one stipulates that $\Dcal(z)^T\eta(u)=\eta(\Dcal(z)u)$, for all $u\in \bw^r\hskip-2pt V$. By 
\cite[Proposition 2.8]{pluckercone} $\Dcal(z)^T$ is a HS--derivation of $\bwV^*$.

\section{Recap on Schubert Derivations}\label{sec:sec3}

\claim{}  Recall that there is a natural representation $\delta:gl(V)\sra \End(\bwV)$ 
making any $\phi\in gl(V)$ into an (even) derivation $\delta(\phi)$ of $\bwV$. In other 
words $\delta(\phi)$ is the unique $\QQ$-vector space endomorphism of $\bwV$ such that
$$
\delta(\phi)(v\w w)=\delta(\phi)v\w w+v\w \delta(\phi)w
$$
for all $v,w\in\bwV$ together with the initial condition  $\delta(\phi)u=\phi(u)$ for all $u\in V=\bw^1V$.
An easy check shows that
$$
\Dcal^\phi(z)=\exp\left(\sum_{i\geq 1}{1\over i}\delta(\phi^i)z^i\right),
$$
is the unique HS derivation on $\bwV$ such that
$\Dcal^\phi(z)_{|_V}=\sum_{i\geq 0}\phi^iz^i$.

\smallskip
Let now $\sigma_1:V\sra V$ be such that $\sigma_1b_j=b_{j+1}$  and $\sigma_{-1}:V\sra V
$ such that $\sigma_{-1}b_j=b_{j-1}$, where by convention we put $b_{k}=0$ if $k<0$.
\claim{\bf Definition.}
{\em The  {\em Schubert derivations} on $\bwV$ are the HS-derivations $\sigma_+(z):\bwV
\sra \bwVz$ and  $\sigma_-(z):\bwV\sra \bwV[z^{-1}]$ defined by
\begin{eqnarray}
\sigma_+(z)&=&\sum_{i\geq 0}\sigma_iz^i:=\exp\left(\sum_{i\geq 1}{1\over i}
\delta(\sigma_1^i)z^i\right),\\  \sigma_-(z)&=&\sum_{i\geq 0}(-1)^i\sigma_{-i}z^{-i}:=
\exp\left(\sum_{i\geq 1}{1\over i}\delta({\sigma_{-1}}\hskip-10pt^i\,\,\,)z^{-i}\right),
\end{eqnarray}
and their inverses in $\End_\QQ(\bwV)\llb z\rrb$ and $\End_\QQ(\bwV)[z^{-1}]$ 
respectively:
\begin{eqnarray}
\ovsig_+(z)&=&\sum_{i\geq 0}(-1)^i\ovsig_iz^i:=\exp\left(-\sum_{i\geq 1}{1\over i}
\delta(\sigma_1^i)z^i\right),\\  \ovsig_-(z)&=&\sum_{i\geq 0}(-1)^i\ovsig_{-i}z^{-i}:=
\exp\left(-\sum_{i\geq 1}{1\over i}\delta({\sigma_{-1}}\hskip-10pt^i\,\,\,)z^{-i}
\right).
\end{eqnarray}
In particular:
$$
\ovsig_{\pm}(z)u=u-\sigma_{\pm 1}u\cdot z^{\pm 1}, \qquad \forall u\in V=\bw^1V.
$$

\medskip

\claim{} It might be useful for the interested reader but not familiar with the main reference of the subject seeing explicitly the action of $\sigma_{\pm}(z)$ and $\ovsig_\pm(z)$ on basis elements.
One has, for all $j\geq 0$
\be
\sigma_+(z)b_j=\sum_{i\geq 0}b_{j+i}z^i \qquad\mathrm{and}\qquad \ovsig_+(z)b_j=b_j-b_{j+1}z,\label{eq0:s+bi}
\ee
\be
\sigma_-(z)b_j=\sum_{i\geq 0}{b_{j-i}\over z^i}, \qquad\mathrm{and}\qquad  \ovsig_-(z)b_j=b_j-{b_{j-1}\over z},\label{eq0:s-bi}
\ee
putting $b_i=0$ for $i<0$. 
\claim{}\label{sec3:mosst} We exploit the  Schubert derivation  $\ovsig_+(z)$ or, equivalently, its 
inverse  $\sigma_+(z)$,  to endow  $\bw^rV$ with a $B_r$-module structure, by declaring 
that
$
e_iu=\ovsig_iu$ or, equivalently, $h_iu=\sigma_iu$,  $\forall u\in \bw^r V
$.
In particular, for all $u\in\bw^rV$:
$$
\ovsig_+(z)u={E_r(z)}\cdot u\qquad\mathrm{and}\qquad \sigma_+(z)u:={1\over E_r(z)}u,
\qquad \forall u\in\bw^rV .
$$
The fact that such a product structure is compatible with the natural vector space 
isomorphism $B_r\sra \bw^r V$ given by \eqref{eq1:fbfc} is a consequence of
\bclm{\bf Proposition.}\label{prop:giamb} {\em Giambelli's formula for the Schubert 
derivation $\sigma_+(z)$ holds:
\be
\wb^r_\blamb=\Delta_\blamb(\sigma_+(z)):=\det(\sigma_{\lambda_j-j+i})_{1\leq i,j\leq r}\wb^r_0.
\label{eq:gmbsh}
\ee
Hence $\bw^rV$ is a free $B_r$-module of rank $1$ generated by $\wb^r_0$.
}
\eclm
\proof Formula~\eqref{eq:gmbsh}  may be inferred as a particular case of the general 
determinantal formula for the exterior power of a polynomial ring due to Laksov and 
Thorup as in  \cite[Main Theorem 0.1]{LakTh01}. It follows that  $\wb^r_\blamb=\Delta_
\blamb(\sigma_+(z))\wb^r_0=\Delta_\blamb(H_r)\wb^r_0$, proving the second part of the 
claim.\qed

\claim{}\label{sec:36} By virtue of \ref{prop:giamb},  the map $B_r\sra \bw^rV$ defined by $\Delta_
\blamb(H_r)\mapsto \Delta_\blamb(H_r)\wb^r_0$ is well defined and is an isomorphism as it maps the basis $(\Delta_\blamb(H_r))_{\blamb\in\Pcal_r}$ of $B_r$ to the basis $(\wb^r_\blamb)_{\blamb\in\Pcal_r}$ of $\bw^rV$.
\claim{} The fact that $\bw^rV$ is a free $B_r$-module of rank $1$ generated by $\wb^r_0$, as prescribed by equality \eqref{eq:gmbsh}, shows that the Schubert derivations  $\sigma_-(z), \ovsig_-(z)$ induces maps $B_r\sra B_r[z^{-1}]$ that, abusing notation, will be denoted in the same way. Their action on a basis element $\Delta_\blamb(H_r)$ of $B_r$ is defined through its action on $\bw^rV$:
\begin{eqnarray}
(\ovsig_-(z)\Delta_\blamb(H_r))\wb^r_0&=&\ovsig_-(z)\wb^r_\blamb,\label{eq3:ovsig-b}\\ 
\cr (\sigma_-(z)\Delta_\blamb(H_r))\wb^r_0&=&\sigma_-(z)\wb^r_\blamb.\label{eq3:sig-b}
\end{eqnarray}

Denote by $\ovsig_-(z)H_r$ (respectively $\sigma_-(z)H_r$)  the sequence  $(\ovsig_-(z)h_j)_{j\in\ZZ}$ (respectively $(\sigma_-(z)h_j)_{j\in\ZZ}$). Then the  following statement gives a practical way to evaluate the image of $\Delta_\blamb(H_r)$ through the maps $\ovsig_-(z)$ and $\sigma_-(z)$ defined by \eqref{eq3:sig-b} and \eqref{eq3:ovsig-b}.
\claim{\bf Proposition.} \label{prop3:prop35} {\em 
For all $r\geq 0$ and all $\blamb\in\Pcal_r$
\be 
\sigma_-(z)h_j=\sum_{i\geq 0}{h_{j-i}\over z^i}\qquad \mathrm{and}\qquad \ovsig_-(z)h_j=h_j-{h_{j-1}\over z}.
\ee
Moreover:
\be
\sigma_-(z)\Delta_\blamb(H_r)=\Delta_\blamb(\sigma_-(z)H_r)\qquad\mathrm{and}\qquad \ovsig_-(z)\Delta_\blamb(H_r)=\Delta_\blamb(\ovsig_-(z)H_r).\label{eq1:comm1}
\ee
}

\proof  See \cite[Theorem 5.7]{pluckercone}, by exploiting the Laksov \& Thorup determinantal formula as in \cite[Main Theorem 0.1]{LakTh01}.\qed

\claim{\bf Remark.}	It is important to notice that \eqref{eq1:comm1} only holds if $\ell(\blamb)\leq r$. For example
$$
\Delta_{(1,1)}(\ovsig_-(z)H_1)=
\left|\bmtx h_1-\displaystyle{1\over z}&1\cr\cr
h_2-\displaystyle{h_1\over z}& h_1-\displaystyle{1\over z}\emtx\right|=-{h_1\over z}+{1\over z^2}\neq 0=\ovsig_-(z)\Delta_{(1,1)}(H_1).
$$

\vspace{-5pt}

%

%\section{Deformation of Schur determinants}
% 
%
%\claim{\bf Notation.}  Let $\bfw_k:=(w_1,\ldots, w_k)$ be an ordered $k$-tuple of indeterminates over $\QQ$. For all $s\geq 1$,  let $H_s:=(h_j)_{j\in\ZZ}$ be the sequence defined by the inverse of the polynomial $E_s(z)$ as in Section. For all partition $\blamb\in\Pcal_r$, and all $0\leq k\leq r$,  the following notation will be adopted:
%%	$$
%%	\Delta_\blamb(\bfw_k^{-1},H_s)=\Delta_\blamb((w_1^{-j}),\ldots,w_k^{(-j)}, \underbrace{H_s,\ldots, H_s}_{r-k\,\,\mathrm{times}}))
%%	$$
%%	will be adopted. 
%%	More explictly:
%$$
%\Delta_\blamb(\bfw^{-1}_k,H_s):=\begin{vmatrix}
%w_1^{-\lambda_1} & w_1^{1-\lambda_2} &\ldots & w_1^{r-1-\lambda_r} \\
%\vdots  & \vdots & \ddots & \vdots \\
%w_k^{-\lambda_1-k+1} & w_k^{-\lambda_2-k+2} &\ldots & w_k^{-\lambda_r+k-r} \\
%& & & \\
%h_{\lambda_1+k} & h_{\lambda_2+k-1} &\ldots &  h_{\lambda_r+k-r+1} \\
%\vdots  & \vdots & \ddots & \vdots \\
%h_{\lambda_1+r-1} & h_{\lambda_2+r-2} &\ldots &  h_{\lambda_r}  \\
%\end{vmatrix}.
%$$
%
%
%\claim{\bf Example} In case $k=0$, one has $\Delta_\blamb(H_s^{(r)})=\det(h_{\lambda_j-j+i})_{1\leq i,j\leq r}$ is the usual Schur determinant. In case $k=r$:
%$$
%\Delta_\blamb(\bfw_r^{-1})=\prod_{j=1}^rw_j^{r-j}s_\blamb(\bfw_k^{-1})\Delta_0(\bfw_k)
%$$
%where 
%$$
%\Delta_0(\bfw_k)=\prod_{i<j} (w_i^{-1}-w_j^{-1})
%$$
%is nothing but the Vandermonde determinant $\prod_{i<j}(w_i^{-1}-w_j^{-1})$ and $s_\blamb(\bfw^{-1}_k)$ is the Schur polynomial associated to the partition $\blamb$ and the  variables $(w_1^{-1},\ldots,w_k^{-1})$.
%
%
%%

\section{Commutation rules}\label{sec:sec4}
\claim{}\label{sec4:sec41} Let $k\geq 1$ and let $\bfz_k:=(z_1,\ldots,z_k)$ be one ordered $k$-tuple  of formal variables. By $\bfz_k^{-1}$ we shall mean the $k$-tuple of the formal inverses $(z_1^{-1},\ldots, z_k^{-1})$.  Define maps $\sigma_\pm(\bfz_k),\ovsig_\pm(\bfz_k):\bwV\sra\bwV\llb\bfz_k,\bfz_k^{-1}]$ respectively by
\be
\sigma_\pm(\bfz_k):=\sigma_\pm(z_1)\cdots\sigma_\pm(z_k)\qquad \mathrm{and}\qquad \ovsig_\pm(\bfz_k):=\ovsig_\pm(z_1)\cdots\ovsig_\pm(z_k).\label{eq4:szk}
\ee
The maps occurring in formulas \eqref{eq4:szk} are multivariate HS derivations on $\bwV$, in the sense that, for instance, $\sigma_+(\bfz_k)(u\w v)=\sigma_+(\bfz_k)u\w \sigma_+(\bfz_k)v$, as it is easy to check and adopting the linear extension of the Schubert derivation to polynomial coefficients as recalled in Definition \ref{def:HSDr}. The same holds verbatim for  $\sigma_-(\bfz_k)$ and $\ovsig_\pm(\bfz_k)$. It is an important point that the multivariate HS derivations  in \eqref{eq4:szk} are symmetric in the formal variables $z_i$ and $w_i$. This is a consequence of the first of the commutation rules of product of Schubert derivations we shall list in this section because needed in the sequel.
\bclm{\bf Proposition.}\label{commu} {\em Let $z,w$ be arbitrary formal variables. The equalities 
\medskip
\begin{eqnarray}
\ovsig_\pm(z)\ovsig_\pm(w)&=&\ovsig_\pm(w)\ovsig_\pm(z)\label{eq:fulcom1},\\ \cr
\sigma_\pm(z)\sigma_\pm(w)&=&\sigma_\pm(w)\sigma_\pm(z)\label{eq:fulcom2}.
\end{eqnarray}
hold  in $\End_\QQ(\bwV)[[z^{\pm},w^\pm]]$.
}
\eclm
\proof Equalities \eqref{eq:fulcom1} and \eqref{eq:fulcom2} hold if and only if  are obvious consequences of the fact if $i,j\geq 0$ then $\sigma_{\pm i}$ and $\sigma_{\pm j}$
are pairwise commuting. It is sufficient then to show that they commute when restricted to $V$, because if they do, then 
$$
\sigma_\pm(z)\sigma_\pm(w)\wb^r_\blamb=\Big[\sigma_\pm(z)\sigma_\pm(w)\bfb\Big]^r_\blamb=\Big[\sigma_\pm(w)\sigma_\pm(z)\bfb\Big]^r_\blamb,
$$ with obvious meaning of the notation. But $\sigma_{\pm i}\sigma_{\pm j}u=\sigma_{\pm 1}^{i+j}u=\sigma_{\pm j}\sigma_{\pm i}u$ for all $u\in V$ and then the claim follows.\qed

In order to give a compact expression of the $gl(\bw^kV)$-module structure of $B_r$, we shall need to introduce a generalisation of the classical vertex operators arising in the context of the so--called boson-fermion correspondence, like in e.g. \cite{KacRaRoz}, which we look at as a generalisation of the isomorphism $B_r\sra \bw^rV$ recalled in Section~\ref{sec2:22}, reaffirmed and refined in Section~\ref{sec:36}.
\claim{\bf Definition.}\label{def4:vertex} {\em By {\em vertex operators} on $\bwV$ we  mean the $\QQ[\bfz_k,\bfz_k^{-1}]$ --linear maps \linebreak$\Gamma(\bfz_k),\Gamma^*(\bfz_k):\bwV\sra (\bwV)\llb \bfz_k, \bfz_k^{-1}]$ of degree $1$ and $-1$, with respect to the exterior algebra graduation, given by:
\begin{eqnarray}
\Gamma(\bfz_k)\wb^r_\blamb&=&\sigma_+(\bfz_k)\ovsig_-(\bfz_k)\wb^{r+k}_\blamb,\label{eq4:gamma}\\ \cr
\Gamma^*(\bfz_k)\wb^r_\blamb&=&\Big(\ovsig_+(\bfz_k)\Delta_\blamb(\sigma_-(\bfz_k)H_{r-k})\Big)\wb^{r-k}_0.\label{eq4:gamma*}
\end{eqnarray}
}
\noindent
\claim{} Proposition \ref{commu} guarantees that the vertex operators $\Gamma(\bfz_k)$ and $\Gamma^*(\bfz_k)$ are symmetric in the formal variables $(z_1,\ldots,z_k)$. They will be studied in a more detailed way in Section \ref{sec:sec5} and \ref{sec:sec6}, exploiting further  commutation relations, for which we need the preliminary work exposed below. As a matter of fact,  we notice that the commutativity of the product of Schubert derivations is granted only if they are of the same kind (both subscript ``$+$'' or both subscript ''$-$''). In general, for $i,j>0$, $\sigma_i$ and $\sigma_{-j}$ do not commute, because $\sigma_{-j}$ is locally nilpotent. The simplest example is: $\sigma_{-1}\sigma_1b_0=b_0\neq 0=\sigma_1\sigma_{-1}b_0$. The general pattern is that commutativity only holds up to the multiplication by a rational function. A first instance of  non trivial commutation rule, needed in the sequel, is provided by the following:

\bclm{\bf Proposition.} \label{prop4:propo44}{\em
\begin{enumerate}
\item[i)] If $\blamb\in \Pcal_r\setminus\Pcal_{r-1}$ {\em (i.e. $\ell(\blamb)=r$)},  then $\ovsig_-(w)$ commutes with both $\sigma_+(z)$ and $\ovsig_+(z)$, i.e.
\be
\ovsig_-(w)\sigma_+(z)=\sigma_+(z)\ovsig_-(w),\label{eq4:com03}
\ee
and
\be
\ovsig_-(w)\ovsig_+(z)=\ovsig_+(z)\ovsig_-(w).\label{eq4:com04}
\ee
\item[ii)] if $\blamb\in\Pcal_{r-1}$ {\em (}i.e. $\wb^r_\blamb=\wb^{r-1}_{\blamb+(1^{r-1})}\w b_0${\em )}{\em :}
\begin{eqnarray}
\ovsig_-(w)\sigma_+(z)\wb^{r}_\blamb&=&\left(1-\displaystyle{z\over w}\right)\sigma_+(z)\ovsig_-(w)\wb^{r}_\blamb\label{eq:com04}.
\end{eqnarray}
\end{enumerate}
}
\eclm
\proof As a matter of i), we observe that $\ovsig_-(w)\sigma_+(z)b_\lambda=\sigma_+(z)\ovsig_-(z)b_\lambda$ if $\lambda>0$. Indeed
\begin{center}
	\begin{tabular}{rllrr}
	$\ovsig_-(w)\sigma_+(z)b_\lambda$&$=$&$\ovsig_-(w)\left(\displaystyle{\sum_{j\geq 0}}{b_{\lambda+j}z^j}\right)$&&\hskip 60pt(Definition of $\sigma_+(z)b_\lambda$)\cr\cr
	&$=$&$\displaystyle{\sum_{j\geq 0}}\left(b_{\lambda+j}-\displaystyle{b_{\lambda+j-1}\over w}\right)z^j$&&\hskip 50pt(Definition of $\ovsig_-(w)$)\cr\cr
	&$=$&$\sigma_+(z)\ovsig_-(w)b_\lambda$.	
	\end{tabular}
\end{center}

%\begin{eqnarray*}
%\ovsig_-(w)\sigma_+(z)b_\lambda&=&\ovsig_-(w)\left(\sum_{j\geq 0}{b_{\lambda+j}z^j}\right)\cr\cr
%&=&\sum_{j\geq 0}\left(b_{\lambda+j}-{b_{\lambda-j-1}\over w}\right)z^j\cr\cr
%&=&\sigma_+(z)\ovsig_-(w)b_\lambda.
%\end{eqnarray*}

Similarly
$$
\ovsig_-(w)\ovsig_+(z)b_\lambda=\ovsig_+(z)\ovsig_-(w)b_\lambda,
$$
as a direct straightforward computation shows. Therefore, under the hypothesis $\ell(\blamb)=r$:
\begin{eqnarray*}
\ovsig_-(w)\sigma_+(z)\wb^r_\blamb&=&\ovsig_-(w)\sigma_+(z)b_{r-1+\lambda_1}\w\cdots\w \ovsig_-(w)\sigma_+(z)b_{\lambda_r}\cr\cr
&=&\sigma_+(z)\ovsig_-(w)b_{r-1+\lambda_1}\w\cdots\w \sigma_+(z)\ovsig_-(w)b_{\lambda_r}=\sigma_+(z)\ovsig_-(w)\wb^r_\blamb,
\end{eqnarray*}
and the same can be argued for the commutation of $\ovsig_+(z)$ and $\ovsig_-(w)$.

\smallskip
To prove ii), equality \eqref{eq:com04}, one observes that
\begin{center}
\begin{tabular}{rllrr}
$\ovsig_-(w)\sigma_+(z)b_0$&$=$&$\ovsig_-(w)\sum_{j\geq 0}b_jz^j$&&\hskip 40pt(Definition of $\sigma_+(z)b_0$)\cr\cr
&$=$&$b_0+ \displaystyle{\sum_{j\geq 1}}\left(b_{j}-\displaystyle{b_{j-1}\over w}\right)z^j$\cr\cr
&$=$&$b_0+\displaystyle{\sum_{j\geq 1}}b_jz^j-\displaystyle{z\over w}\sum_{j\geq 0}b_jz^j$\cr\cr
&$=$&$\left(1-\displaystyle{z\over w}\right)\sigma_+(z)b_0$
\end{tabular}
\end{center}
\be
\hskip-108pt =\,\,\,\,\left(1-\displaystyle{z\over w}\right)\sigma_+(z)\ovsig_-(w)b_0
\label{Prop:45}
\ee
%\begin{center}
%\be
%\begin{tabular}{rllrr}
%$=$&$\left(1-\displaystyle{z\over w}\right)\sigma_+(z)\ovsig_-(w)b_0$,&\qquad\qquad\qquad\qquad\qquad\qquad&\label{Prop:45}	\end{tabular}
%\ee	
%\end{center}
\smallskip
because, in general, $\ovsig_-(w)$ acts on $\wb^r_0$ as the identity.
So,  if $\ell(\blamb)<r$ (i.e. $\lambda_r=0$) one obtains:
\begin{center}
	\begin{tabular}{rclrr}
\hskip-3pt $\ovsig_-(w)\sigma_+(z)\wb^r_\blamb$&\hskip-10pt$=$&$\hskip-10pt\ovsig_-(w)\sigma_+(z)\left(\wb^{r-1}_{\blamb+(1^{r-1})}\w b_0\right)$&&\hskip-8pt (Definition\cr
&&&& of $\wb^r_\blamb$)\cr\cr
&\hskip-9pt$=$&\hskip-10pt $\ovsig_-(w)\sigma_+(z)\wb^{r-1}_{\blamb+(1^{r-1})}\w \ovsig_-(w)\sigma_+(z)b_0$&&($\ovsig_-(w)\sigma_+(z)$ is\cr 
&&&& a HS derivation)\cr\cr
&\hskip-9pt $=$&\hskip-10pt$\sigma_+(z)\ovsig_-(w)\wb^{r-1}_{\blamb+(1^{r-1})}\w\left(1-\displaystyle{z\over w}\right)\sigma_+(z)\ovsig_-(w)b_0$&&Commu--\cr &&&&tation \eqref{Prop:45})\cr\cr
&\hskip-9pt$=$&\hskip-10pt$\left(1-\displaystyle{z\over w}\right)\sigma_+(z)\ovsig_-(w)\wb^r_\blamb$.&&\qed
\end{tabular}
\end{center}

The rest of this section will be devoted to prove the following 
\bclm{\bf Theorem.}\label{thm:thm44} {\em For all $u\in\bw^rV$, the following commutation rule holds:
\be
\beta_0\lrcorner \sigma_-(w)\ovsig_+(z)u=\left(1-{z\over w}\right)\ovsig_+(z)(\beta_0\lrcorner \sigma_-(w)u).
\ee
}
\eclm
\claim{} \label{sec:sec46} Some preparation is needed to prove Theorem \ref{thm:thm44}. First we introduce a piece of useful notation. 
If $\beta\in V^*$, the contraction of the decomposable tensor $u_1\w\cdots\w u_r$ of $\bw^rV$ against $\beta$ may be represented via the following diagram:
\be
\begin{vmatrix}\beta\lrcorner u_1& \beta\lrcorner u_2&\cdots&\beta\lrcorner u_r\cr\cr
u_1&u_2&\ldots&u_r\end{vmatrix}=\begin{vmatrix}\beta(u_1)& \beta(u_2)&\cdots&\beta(u_r)\cr\cr
u_1&u_2&\ldots&u_r\end{vmatrix},\label{eq4:contdiag}
\ee
which  precisely means that $(-1)^{j+1} \beta\lrcorner u_j$ is the coefficient of the vector of $\bw^{r-1}V$  obtained  by removing the wedge factor $u_j$ from   $u_1\w u_2\w\cdots\w u_r$.
For example 
$$
\begin{vmatrix} \beta\lrcorner u_1 &\beta\lrcorner u_2&\beta\lrcorner u_3\cr\cr
u_1&u_2&u_3
\end{vmatrix}=\beta(u_1)\cdot u_2\w u_3-\beta(u_2) u_1\w u_3+\beta(u_3) u_1\w u_2,
$$
which is precisely the expanded expression of the contraction $\beta\lrcorner (u_1\w u_2\w u_3)$.

\claim{} Recall   the generating 
function $\bbeta(w^{-1}):=\displaystyle{\sum_{j\geq 0}}\beta_jw^{-j}$ introduced in formula \eqref{eq2:generating}.  Since $\beta_j\lrcorner b_i=\beta_j(b_i)=\delta_{ij}$, it clearly follows that
$\bbeta(w^{-1})\lrcorner b_j=\displaystyle{\sum}_{i\geq 0}\beta_i(b_j)w^{-i}=w^{-j}$. Therefore:
\be 
\bbeta(w^{-1})\lrcorner\wb^r_\blamb=\begin{vmatrix} w^{-r+1-\lambda_1}&w^{-r+2-\lambda_2}&\cdots&w^{-\lambda_r}\cr\cr
b_{r-1+\lambda_1}&b_{r-2+\lambda_2}&\cdots&b_{\lambda_r}
\end{vmatrix}.
\ee

\bclm{\bf Proposition.} \label{prop:prop48}{\em For all
	$u\in\bw^rV$:
\be	
\bbeta(w^{-1})\lrcorner \ovsig_+ (z)u=\left(1-{z\over w}\right)\ovsig_+(z)(\bbeta(w^{-1})\lrcorner u).
\ee
}
\eclm
\proof	Since each $u\in\bw^rV$ is a finite linear combination of $\wb^r_\blamb$, it is no harm to assume $u=\wb^r_\blamb$. Then we start to notice that
\begin{center}
	\begin{tabular}{rclr}
$\bbeta(w^{-1})\lrcorner \ovsig_+(z)b_j$&$=$&$\bbeta(w^{-1})\lrcorner (b_j-b_{j+1}z)$&\hskip 70pt (Definition of $\ovsig_+(z)$)\cr\cr
&$=$&$\bbeta(w^{-1})\lrcorner b_j- \bbeta(w^{-1})\lrcorner b_{j+1}z$& \hskip 30pt (Action of $\lrcorner$)
\end{tabular}
\end{center}
\begin{center}
	\be
	\begin{tabular}{rllrr}
	&$=$&$\displaystyle{1\over w^{j}}-\displaystyle{z\over w^{j+1}}=\displaystyle{1\over w^j}\left(1-\displaystyle{z\over w}\right)$.\qquad\qquad\qquad\qquad\qquad\label{eq4:auto}
\end{tabular}
\ee
\end{center}
By using the expression of a contraction via diagram \eqref{eq4:contdiag}, one has:
\begin{eqnarray*}
&&\bbeta(w^{-1})\lrcorner \ovsig_+ (z)\wb^r _{\blamb}\cr\cr
&=&\begin{vmatrix}\bbeta(w^{-1})\lrcorner\ovsig_+(z)b_{r-1+\lambda_1}&\bbeta(w^{-1})\lrcorner\ovsig_+(z)b_{r-2+\lambda_2}&\cdots&\bbeta(w^{-1})\lrcorner\ovsig_+(z)b_{\lambda_r}\cr\cr\cr
\ovsig_+(z)b_{r-1+\lambda_1}&\ovsig_+(z)b_{r-2+\lambda_2}&\cdots&\ovsig_+(z)b_{\lambda_r}\end{vmatrix}.
\end{eqnarray*}
\smallskip
which by \eqref{eq4:auto} is equal to:
\smallskip
\begin{eqnarray}
&=&\begin{vmatrix}\left(1-\displaystyle{z\over w}\right)\displaystyle{1\over w^{r-1+\lambda_1}}&\left(1-\displaystyle{z\over w}\right)\displaystyle{1\over w^{r-2+\lambda_2}}&\cdots&\left(1-\displaystyle{z\over w}\right)\displaystyle{1\over w^{\lambda_r}}\cr\cr\cr
\ovsig_+(z)b_{r-1+\lambda_1}&\ovsig_+(z)b_{r-2+\lambda_2}&\cdots&\ovsig_+(z)b_{\lambda_r}\end{vmatrix}\cr\cr\cr
&=&\left(1-\displaystyle{z\over w}\right)\begin{vmatrix}\displaystyle{1\over w^{r-1+\lambda_1}}&\displaystyle{1\over w^{r-2+\lambda_2}}&\cdots&\displaystyle{1\over w^{\lambda_r}}\cr\cr
\ovsig_+(z)b_{r-1+\lambda_1}&\ovsig_+(z)b_{r-2+\lambda_2}&\cdots&\ovsig_+(z)b_{\lambda_r}\end{vmatrix}.\label{eq4:ourdet}
\end{eqnarray}
Since the determinant occurring in \eqref{eq4:ourdet} is a linear combination of 
$\big[\ovsig_+(z)\bfb\big]^{r-1}_{\blamb^{(j)}}=\ovsig_+(z)\wb^{r-1}_{\blamb^{(j)}}$ (because $\ovsig_+(z)$ is a HS derivation), where we denoted by $\blamb^{(j)}$ the partition of lenght at most $r-1$ obtained by omitting the $j$-th part, it follows that the action of $\ovsig_+(z)$ can be factorized from the bottom row of \eqref{eq4:ourdet}, giving
$$
\left(1-\displaystyle{z\over w}\right)\ovsig_+(z)\begin{vmatrix}\displaystyle{1\over w^{r-1+\lambda_1}}&\displaystyle{1\over w^{r-2+\lambda_2}}&\cdots&\displaystyle{1\over w^{\lambda_r}}\cr\cr
b_{r-1+\lambda_1}&b_{r-2+\lambda_2}&\cdots&b_{\lambda_r}\end{vmatrix}=\left(1-\displaystyle{z\over w}\right)\ovsig_+(z)(\bbeta(w^{-1})\lrcorner \wb^r_\blamb),
$$
which ends the proof of the Proposition. \qed
\bclm{\bf Lemma.} \label{lem4:35}{\em  
For all $u\in\bw^rV$, 
\be
\bbeta(w^{-1})\lrcorner u=\ovsig_-(w)(\beta_0\lrcorner \sigma_-(w)u).\label{eq4:gsal}
\ee
}
\eclm
\proof	The proof is basically contained in \cite[Proposition 4.3]{pluckercone} but, because some mild deformity in the  notation, we prefer to repeat it here.
Recall the definition of transpose of a HS derivation on $\bw V$. We observe that $\bbeta(w^{-1})=\sigma_-(z)^T\beta_0$. Then, for all $\eta\in\bw^{r-1}V$, 
\begin{center}
	\begin{tabular}{rclr}
$\eta(\bbeta(w^{-1})\lrcorner u)$&$=$&$(\bbeta(w^{-1})\w \eta)(u)$&(Definition \ref{pairing} of contraction)\cr\cr
&$=$&$(\sigma_-(w)^T\beta_0\w\eta)(u)$&(by the above observation)\cr\cr
&$=$&$\sigma_-(w)^T(\beta_0\w\ovsig_-(w)^T\eta)(u)$&(Integration by parts)\cr\cr %$\mathrm{(integration\,\,by\,\, parts)}$\cr\cr
&$=$&$\beta_0\w\ovsig_-(w)^T\eta(\sigma_-(w)u)$&(Definition of transpose of $\ovsig_-(w)^T$)\cr\cr
&$=$&$\ovsig_-(w)^T\eta(\beta_0\lrcorner \sigma_-(w)u)$&(Again \ref{pairing})\cr\cr
&$=$&$\eta(\ovsig_-(w)(\beta_0\lrcorner \sigma_-(w)u))$&(Definition of transpose).
\end{tabular}
\end{center}
The last equality proves \eqref{eq4:gsal}, due to the arbitrary choice of $\eta\in \bw^{r-1}V^*\cong(\bw^{r-1}V)^*$.\qed

\bclm{\bf Lemma.} \label{lem4:411}{\em Taking $\ovsig_+(z)$ commutes with takin contraction against $\beta_0$, i.e. for all $\blamb\in\Pcal_r$
\be
\beta_0\lrcorner \ovsig_+(z)\wb^r_\blamb=\ovsig_+(z)(\beta_0\lrcorner \wb^r_\blamb).\label{eq4:s+cb}
\ee
}
\eclm
\proof There are two cases. If $\ell(\blamb)=r$, both members of \eqref{eq4:s+cb} vanish.
If $\ell(\blamb)\leq r-1$, then
$\beta_0\lrcorner \ovsig_+(z)\wb^r_\blamb=(-1)^{r-1}\ovsig_+(z)\wb^{r-1}_{\blamb+(1^{(r-1)})}=\ovsig_+(z)(\beta_0\lrcorner \wb^r_\blamb)$.
\qed

\smallskip
We are now in position to  provide the

\claim{\bf Proof of Theorem \ref{thm:thm44}.} We have:

\begin{center}
\begin{tabular}{rllrr}
$\beta_0\lrcorner \sigma_-(w)\ovsig_+(z)u$&$\hskip-7pt=$&$\hskip-11pt\sigma_-(w)\ovsig_-(w)(\beta_0\lrcorner \sigma_-(w)\ovsig_+(z)u)$&&($\sigma_-(w)\ovsig_-(w)=1$)\cr\cr
&$\hskip-7pt=$&$\hskip-11pt\sigma_-(w)(\bbeta(w^{-1})\lrcorner\ovsig_+(z)u)$&&(Lemma \ref{lem4:35})\cr\cr
&$\hskip-7pt=$&$\hskip-11pt\left(1-\displaystyle{z\over w}\right)\sigma_-(w)\ovsig_+(z)(\bbeta(w^{-1})\w u)$&&(Proposition \ref{prop:prop48})\cr\cr
&$\hskip-7pt=$&$\hskip-11pt\left(1-\displaystyle{z\over w}\right)\sigma_-(w)\ovsig_+(z)\ovsig_-(w)(\beta_0\lrcorner \sigma_-(w)u)$&&(again Lemma \ref{lem4:35})\cr\cr
&$\hskip-7pt=$&$\hskip-11pt\left(1-\displaystyle{z\over w}\right)\sigma_-(w)\ovsig_-(w)\ovsig_+(z)(\beta_0\lrcorner \sigma_-(w)u)$&&(formula \eqref{eq4:com04} \cr
&&&&of Proposition \ref{prop4:propo44})\cr\cr
&$\hskip-7pt=$&$\hskip-13pt\left(1-\displaystyle{z\over w}\right)\ovsig_+(z)(\beta_0\lrcorner \sigma_-(w)u)$&&($\sigma_-(w)\ovsig_-(w)=1$).
\end{tabular}
\end{center}
and Theorem \ref{thm:thm44}  is thence proven. \qed

%%\section{Vertex operators on exterior algebras}
%%
%%Notation as in Section \ref{sec4:sec41} and recall Proposition \ref{prop3:prop35}, notably \eqref{eq1:comm1}.
%%
%%
%%\bclm{\bf Proposition.}
%%
%%\eclm
%%
%% 
%% \claim{\bf Proposition.} \label{prop:prop52}{\em
%% One has 
%%\be
%%\Gamma(z)\wb^r_\blamb=z^{-r}\sigma_+(z)b_0\w \wb^r_\blamb={z^{-r}\over E_r(z)}\ovsig_-(z)\Delta_\blamb(H_{r+1})\wb^{r+1}_0\label{eq:voint}
%%\ee
%%and
%%\be 
%%\Gamma^*(w)\Delta_\blamb(H_r)\wb^r_0=w^{r-1}\bbeta(w^{-1})\lrcorner \wb^r_\blamb=w^{r-1}\Delta_\blamb(\bfw_1^{-1},H_{r-1})\wb^{r-1}_0
%%\ee\label{eq1:vo*int}
%% }
%% 
%% \proof By \cite[Definition 3.5. and Proposition 3.6]{gln}. Notice the use of a different convention  with respect to \cite{gln}. By $\bbeta(\bfw)$ we mean the formal power series $\sum_{j\geq 0}\beta_jw^{-j}$ and not $\sum_{j\geq 0}\beta_jw^{-j-1}$.\qed
% 
%Abusing  notation Formula~\eqref{eq1:vo*int} will be written as
%$$
%\Gamma^*(w)\Delta_\blamb(H_r)=\ovsig_+(w)\Delta_\blamb(\sigma_-(w)H_{r-1})\in B_{r-1}[w^{-1}]
%$$
%In the sequel we  need to consider product of vertex operators of the form $\Gamma(z_1)
%\cdots\Gamma(z_k)$ or $\Gamma^*(w_1)\cdots\Gamma^*(w_k)$. For this reason we  need to state 
%the commutation rules  listed below, which are  in turn related with the locality axioms 
%required in the definition of vertex operators for a vertex algebra. Recall that $\wb^r_\blamb=0$ if $\ell(\blamb)>r$.
\section{The vertex operator $\Gamma(\bfz_k)$}\label{sec:sec5}
The main purpose of this section is to interpret the vertex operator $\Gamma(\bfz_k)$, introduced in Definition \ref{def4:vertex}, formula~\eqref{eq4:gamma}, in terms of wedging operation on the exterior algebra. This generalises \cite[Proposition 4.2]{pluckercone}. This will be done in Theorem \ref{thm:thm61} below and will be used in our main Theorem \ref{thm:mainthm}.  

\bclm{\bf Lemma.} {\em For all $j\geq 0$ and all $k\geq 1$ one has
\begin{eqnarray}
\ovsig_+(\bfz_k)b_j&=&b_j+\sum_{i=1}^k(-1)^ie_i(\bfz_k)b_{j+i}\label{eq1:sim1}\\ \cr
&&\mathrm{and}\cr\cr
\sigma_+(\bfz_k)b_j&=&b_j+\sum_{i\geq 1}^kh_i(\bfz_k)b_{j+i},\label{eq2:sim2}
\end{eqnarray}
where $e_i(\bfz_k)$ and $h_i(\bfz_k)$ are, respectively, the elementary and complete symmetric polynomial of degree $i$ in the indeterminates $\bfz_k:=(z_1,\ldots,z_k)$.
}\eclm
\proof
Formula \eqref{eq1:sim1} is the content of \cite[Lemma 5.7]{SDIWP} to which we refer to. Formula \eqref{eq2:sim2} is a consequence of \eqref{eq1:sim1}, keeping into account that $\sigma_+(\bfz_k)$ and $\ovsig_+(\bfz_k)$ are one the inverse of the other in $\End_\QQ(\bwV)\llb \bfz_k\rrb$.\qed

\bclm{\bf Lemma.}\label{lemma:53} {\em One has:
\begin{eqnarray}
\ovsig_-(\bfz_k)b_{j+k}&=&b_{j+k}+\sum_{i=1}^k(-1)^ie_i(\bfz_k^{-1})b_{j+k-i},\label{eq1:sim-1}
\end{eqnarray}
where $e_i(\bfz_k^{-1})=e_i(z_1^{-1},\ldots,z_k^{-1})$ is the elementary symmetric polynomial of degree $i$ in $(z_1^{-1},\ldots,z_k^{-1})$,
}
\eclm

\proof The  proof works the same as in
  \cite[Lemma 6.7]{SDIWP}. The formula holds for $k=1$ is true, because 
$$
\ovsig_-(z_1)b_{j+1}= b_{j+1}-\displaystyle{b_j \over z_1}.
$$
By induction, suppose that \eqref{eq1:sim-1}  holds for $k-1 \geq 0$. Then it holds for $k$. Indeed
\begin{eqnarray*}
&&\ovsig_-(z_1) \ovsig_-(z_2)\ldots \ovsig_-(z_k)	b_{j+k}\cr\cr
&=&\ovsig_-(z_1)\left[b_{j+k}-e_1 \left({1 \over z_2}, \cdots , {1 \over z_k}\right) b_{j+k-1} + \cdots
+(-1)^{k-1}e_{k-1} \left({1 \over z_2}, \cdots , {1 \over z_k}\right)b_{j+1} \right]\cr\cr
&=& b_{j+k}- {b_{j+k-1}\over z_1} - e_1\left({1 \over z_2}, \cdots , {1 \over z_k}\right)
\left(b_{j+k-1}- {b_{j+k-2}\over z_1} \right)+ \cdots \cr\cr
&+& 
(-1)^{k-1}e_{k-1} \left({1 \over z_2}, \cdots , {1 \over z_k}\right)\left(b_{j+1} -{b_j \over z_1}\right)\qquad\qquad\qquad\qquad \mathrm{(definition\,\,of\,\, }\ovsig_-(z_1))\cr\cr
&=&b_{j+k}+\sum_{i=1}^k(-1)^ie_i\left({1 \over z_1}, \cdots , {1 \over z_k}\right)b_{j+k-i},
\end{eqnarray*} 
as desired.	\qed 

\bclm{\bf Lemma.}\label{lemma34}{\em The following equality holds for all $1\leq i\leq k${\em :}
\be
{e_i(\bfz_k)\over z_1\cdots z_k}=e_{k-i}\left({1\over z_1},\ldots, {1\over z_k}\right). \label{ei}
\ee
}
\eclm
\proof  Recall the following definition of the elementary symmetric polynomials in $k$ indeterminates through generating functions:
\be
\sum_{i=0}^{k}e_i(\bfz_k)t^i = \prod_{i=1}^{k} (1+z_i t).\label{eq5:2sides}
\ee
By dividing both sides of \eqref{eq5:2sides} by $e_k(\bfz_k)=z_1 \ldots z_k$ we get
\be
\sum_{i=0}^{k} {e_i(\bfz_k)\over z_1\cdots z_k}t^i= \prod_{i=1}^{k} \left({1\over z_i}+t\right).\label{eq5:2ssid}
\ee
The claim then follows by comparing the coefficient of $t^i$ in the two sides \eqref{eq5:2ssid}. \qed

\bclm{\bf Lemma.} \label{lem:lemma54}{\em For all $k\geq 1$, $r\geq 0$ and $\blamb\in\Pcal_r$: 
\be
\wb^k_0\w \ovsig_+ (\bfz_k)\wb^{r}_\blamb=e_k(\bfz_k)^r\ovsig_-(\bfz_k)\wb^{r+k}_\blamb.
\label{eq:eii}
\ee
}
\eclm
\proof Equality \eqref{eq:eii} holds for $r=1$:
\smallskip
\begin{center}
	\begin{tabular}{rllrr}
		&&$\wb^k_0\w \ovsig_+(\bfz_k)\wb^{1}_\blamb=
		\wb^k_0 \w \ovsig_+(\bfz_k ) b_\lambda$ & \cr\cr
		&$=$&$\wb^k_0 \w \left( b_\lambda- e_1 (\bfz_k)b_{\lambda+1}+\cdots +
		(-1)^k e_k (\bfz_k) b_{\lambda+k}\right)$&(By \ref {eq1:sim1} )\cr\cr
		&$=$&$\wb^k_0 \w (-1)^k e_k (\bfz_k)\left[ b_{\lambda+k} - \displaystyle{e_{k-1} (\bfz_k) \over e_k (\bfz_k)}b_{\lambda+k-1}+\cdots\right.$\cr
		&&\hskip180pt $	\left. +(-1)^k \displaystyle{1 \over e_k (\bfz_k) }b_{\lambda}\right]$&(By Factorization)\cr\cr
		&$=$&$e_k(\bfz_k)\ovsig_-(\bfz_k) b_{\lambda+k} \w \wb^k_0$&(By lemma \eqref{lemma:53})\cr\cr
		&$=$&$e_k(\bfz_k)\ovsig_-(\bfz_k)( b_{\lambda+k} \w \ovsig_-(\bfz_k) \wb^k_0)$&\big(Integration by parts\big)\cr\cr
		&$=$&$e_k(\bfz_k)\ovsig_-(\bfz_k) \left( b_{\lambda+k} \w b_{k-1} \w \ldots \w b_0 \right)$&($\ovsig_-(\bfz_k) \wb^{k}_0=\wb^{k}_0$)\cr\cr
		&$=$&$ e_k(\bfz_k)\ovsig_-(\bfz_k) \wb^{1+k}_\blamb$.
	\end{tabular}
\end{center}
%
%
%\begin{eqnarray*}
%&&\wb^k_0\w \ovsig_+(\bfz_k)\wb^{1}_\blamb=
%\wb^k_0 \w \ovsig_+(\bfz_k ) b_\lambda & (By \ref {eq1:sim1} )\cr\cr
%&=&\wb^k_0 \w \left( b_\lambda- e_1 (\bfz_k)b_{\lambda+1}+\cdots +
%(-1)^k e_k (\bfz_k) b_{\lambda+k}\right)\cr\cr
%&=&\wb^k_0 \w (-1)^k e_k (\bfz_k)\left[ b_{\lambda+k} - {e_{k-1} (\bfz_k) \over e_k (\bfz_k)}b_{\lambda+k-1}+\cdots +(-1)^k {1 \over e_k (\bfz_k) }b_{\lambda}\right]\cr\cr
%&&&\cr
%&&& (By lemma \ref{lemma34})\cr\cr
%&=&e_k(\bfz_k)\ovsig_-(\bfz_k) b_{\lambda+k} \w \wb^k_0\cr\cr
%&=&e_k(\bfz_k)\ovsig_-(\bfz_k) b_{\lambda+k} \w \ovsig_-(\bfz_k) \wb^k_0\cr\cr
%&=&e_k(\bfz_k)\ovsig_-(\bfz_k) \left( b_{\lambda+k} \w b_k \w \ldots \w b_0 \right)\cr\cr
%&=& e_k(\bfz_k)\ovsig_-(\bfz_k) \wb^{1+k}_\blamb
%\end{eqnarray*}
Therefore the property is true for $r=1$. Assume now \eqref{eq:eii} holds true for $r-1\geq 0$. Then  
\begin{eqnarray*}
\wb^k_0\w \ovsig_+ (\bfz_k)\wb^{r}_\blamb&\hskip-5pt=&\hskip-5pt
\wb^k_0\w \ovsig_+ (\bfz_k) b_{r-1+\lambda_1}\w \ldots \w \ovsig_+ (\bfz_k) b_{\lambda_r}\cr\cr
&\hskip-5pt =& \hskip-5pt \wb^k_0\w (-1)^k e_k(\bfz_k)\ovsig_-(\bfz_k)  b_{r-1+k+\lambda_1}\w \ldots \w (-1)^k e_k(\bfz_k)\ovsig_-(\bfz_k)  b_{\lambda_r +k}\cr\cr
&\hskip-5pt  =&\hskip-5pte_k(\bfz_k)^r \left[ 
\ovsig_-(\bfz_k)  b_{r+k-1+\lambda_1}\w \ldots \w 
\ovsig_-(\bfz_k)  b_{k+\lambda_r }\w \ovsig_-(\bfz_k)\wb^k_0\
\right]\cr\cr
&\hskip-5pt  =& \hskip-5pt e_k(\bfz_k)^r \ovsig_-(\bfz_k)  \left( 
b_{r+k-1+\lambda_1}\w b_{k+\lambda_r } \w b_k \w \ldots \w b_0 \right)\cr\cr
&\hskip-5pt  =&\hskip-5pte_k(\bfz_k)^r\ovsig_-(\bfz_k)\wb^{r+k}_\blamb,
\end{eqnarray*}
as claimed. \qed

\bclm{\bf Theorem.}\label{thm:thm61} {\em For all $u\in \bw^kV\llb \bfw_k,\bfw_k^{-1}]$ we have:
\be 
\sigma_+(z_1,\ldots,z_k)\wb^k_0\w u=\prod_{j=1}^rz_j^r\cdot \Gamma(\bfz_k)u,
\ee
the equality holding  in $\bw^kV\llb\bfz_k,\bfw_k\rrb[\bfw_k^{-1}]$.
}
\eclm
\proof Recall that we consider all the Schubert derivations extended by linearity over rings of formal power series  with rational coefficients. See definition \ref{def:HSDr}. Then our arbitrary $u$ is intended as a possibly infinite linear combination of $\wb^r_\blamb$ with coefficients being polynomials. Then we can assume with no harm that $u=\wb^r_\blamb$, a basis element of $\bw^kV$. Keeping the same notation  as in \ref{sec1:notation}, we first apply integration by parts. Then

\begin{center}
	\begin{tabular}{rclr}
	$\sigma_+(\bfz_k)\wb^k_0 \w 	\wb^r_\blamb$&$=$& 
	$\sigma_+(\bfz_k)\left( \wb^k_0 \w \ovsig_+(\bfz_k)\wb^r_\blamb \right)$\qquad\qquad& (By integration by \cr
	&&&parts \eqref{eq:intp1})\cr\cr
	&$=$&$\sigma_+(\bfz_k)e_k (\bfz_k)^r \ovsig_- (\bfz_k) \wb^{r+k}_\blamb$&(By Lemma \ref{lem:lemma54})\cr\cr
	&$=$&$e_k (\bfz_k)^r \sigma_+(\bfz_k)\ovsig_- (\bfz_k) \wb^{r+k}_\blamb$\cr\cr
	&$=$&$ \prod_{j=1}^kz_j^r\Gamma(\bfz_k)\wb^{r}_\blamb$&(Definition of $\Gamma(\bfz_k)$)
\end{tabular}
 \end{center}

%\begin{eqnarray*}
%\sigma_+(\bfz_k)\wb^k_0 \w 	\wb^r_\blamb&=& 
%\sigma_+(\bfz_k)\left( \wb^k_0 \w \ovsig_+(\bfz_k)\wb^r_\blamb \right)& (By \eqref{eq:eii})\cr\cr
%&=&\sigma_+(\bfz_k)e_k (\bfz_k)^r \ovsig_- (\bfz_k) \wb^{r+k}_\blamb\cr\cr
%&=&e_k (\bfz_k)^r \sigma_+(\bfz_k)\ovsig_- (\bfz_k) \wb^{r+k}_\blamb,\cr\cr
%&=& \prod_{j=1}^rz_j^r\Gamma(\bfz_k)\wb^{r}_\blamb,
%\end{eqnarray*}
as desired. \qed

If $k=1$, and $z=z_1$, one obtains
$$
\sigma_+(z)b_0\w \wb^r_\blamb=z^r\Gamma(z)\wb^{r}_\blamb=z^r\sigma_+(z)\ovsig_-(z)\wb^{r+1}_\blamb,
$$
which is precisely  \cite[Proposition 5.4]{pluckercone} or \cite[Proposition 3.2]{gln}. They shape looks more involved because we use here better notation.

\section{The vertex operator $\Gamma^*(\bfw_k)$}\label{sec:sec6}

In the same vein of Section \ref{sec:sec5}, this section will be devoted to interpret in terms of contraction operations the action of the vertex operator $\Gamma^*(\bfw_k)$ on $\bwV$, The output will be Theorem \ref{thm:thm71}, stated at the end of the section, which will be another building block of the main Theorem \ref{thm:mainthm}. We begin with some preparation.
\bclm{\bf Lemma.}\label{lem:lem62} {\em The following equality holds for  all $r\geq 1$ and all $\blamb\in\Pcal_r${\em :}
\be
\ovsig_{-r+1}(\beta_0\lrcorner \sigma_-(w)\wb^r_\blamb)=\Delta_\blamb(\sigma_-(w)H_{r-1})\wb^{r-1}_0.\label{eq:bw03}
\ee
}
\eclm
\proof This is \cite[Lemma 5.8]{pluckercone}.\qed
%\begin{eqnarray*}
%\end{eqnarray*}
\claim{\bf Lemma.}\label{cor:cor510}
\be
\bbeta(w^{-1})\lrcorner\wb^r_\blamb=w^{-r+1}\ovsig_+(w)\Delta_\blamb(\sigma_-(w)H_{r-1})\wb^{r-1}_0\label{eq:premnthm}
\ee
\proof Invoking Lemma~\ref{lem4:35}, 
\be
\bbeta(w^{-1})\lrcorner\wb^r_\blamb=\ovsig_-(w)(\beta_0\lrcorner \sigma_-(w)\wb^r_\blamb).
\ee
Since $\beta_0\lrcorner \sigma_-(w)\wb^r_\blamb$ is a linear combination of $\wb^{r-1}_\bmu$ with $\ell(\bmu)=r-1$ (i.e. no $b_0$ occurs in the monomial), then by \cite[Proposition 4.3]{pluckercone}
 $$\ovsig_-(w)(\beta_0\lrcorner \sigma_-(w)\wb^r_\blamb)=w^{-r+1}\ovsig_+(w)\ovsig_{-r+1}(\beta_0\lrcorner \sigma_-(w)\wb^r_\blamb).$$
Using Lemma~\ref{lem:lem62} one obtains \eqref{eq:premnthm}.
\qed

\bclm{\bf Theorem.} \label{thm:thm71} {\em The following equality holds:
\be
(\bbeta(w_k^{-1})\w \bbeta(w_{k-1}^{-1})\w\cdots\w \bbeta(w_1 ^{-1}))\lrcorner\wb^r_\blamb={\Delta_0(\bfw_k)\over (w_1\cdots w_k)^{r-1}}\Gamma^*(\bfw_k)\wb^r_\blamb.
\ee
where $\Delta_0(\bfw_k)$ denotes the Vandermonde determinant $\prod_{1\leq i<j\leq k}(w_j-w_i)$.
}
\eclm
\proof
 For $k=1$ the property 
$$
\bbeta(w_1^{-1})\lrcorner\wb^r_\blamb=w_1 ^{r-1}\ovsig_+(w_1)\Delta_\blamb(\sigma_-(w_1)H_{r-1})\wb^{r-1}_0
$$
is just Lemma \ref{cor:cor510}. Arguing by induction, let us assume the claim holding true for $0\leq k-1 \leq r-1$ and let us show it holds for $k$. We have
\begin{eqnarray*}
\bbeta(w_k^{-1})\w \bbeta(w_{k-1}^{-1})\w\cdots\w \bbeta(w_1 ^{-1})\lrcorner\wb^r_\blamb&=&
\bbeta(w_k^{-1}) \lrcorner \left( \bbeta(w_{k-1}^{-1})\w\cdots\w \bbeta(w_1^{-1})\lrcorner\wb^r_\blamb\right).
\end{eqnarray*}	
Using the inductive hypothesis:
\begin{eqnarray*}
&=&\bbeta(w_k^{-1}) \lrcorner
{\Delta_0(\bfw_{k-1})\over (w_1\cdots w_k)^{r-1}}\ovsig_+(\bfw_{k-1})\Delta_\blamb(\sigma_-(\bfw_{k-1})H_{r-k+1}) \wb^{r-k+1}_0.
\end{eqnarray*}
By applying Lemma~\ref{cor:cor510}, one gets:
\begin{eqnarray*}
&=&w_k^{r-k+1}\ovsig_+(w_k)\left[
\beta_0 \lrcorner \sigma_-(w_k)
\ovsig_+(\bfw_{k-1})\Delta_\blamb\big(\sigma_-(\bfw_{k-1})H_{r-k+1}\big) \wb^{r-k+1}_0
\right]\cr\cr
 &&\cdot
{\Delta_0(\bfw_{k-1})\over (w_1\cdots w_{k-1})^{r-1}}
\end{eqnarray*}
Now we use the commutation rules prescribed by Theorem \ref{thm:thm44} and Lemma \ref{lem4:411}:
\begin{eqnarray*}
&=&w_k^{k-r}\ovsig_+(w_k)\left[
\prod_{j=0}^{k-1} \left( 1-{w_j \over w_k}\right)\cdot
\left(\beta_0 \lrcorner 
\ovsig_+(\bfw_{k-1})\sigma_-(\bfw_k) \Delta_\blamb(\sigma_-(\bfw_{k-1})H_{r-k+1}) \wb^{r-k+1}_0 \right)
\right]\cr\cr&& \cdot
{\Delta_0(\bfw_{k-1})\over (w_{1}\cdots w_{k-1})^{r-1}}\cr\cr
&=&{w_k^{k-r}\Delta_0(\bfw_{k-1})\over (w_{1}\cdots w_{k-1})^{r-1}}\prod_{j=0}^{k-1} \left( 1-{w_j \over w_k}\right)\cdot\cr\cr
&&\hskip110pt \cdot\ovsig_+(w_k)\ovsig_+(\bfw_{k-1})\left[
\left(\beta_0 \lrcorner 
\sigma_-(\bfw_k) \Delta_\blamb(\sigma_-(\bfw_{k-1})H_{r-k+1}) \wb^{r-k+1}_0 \right)
\right]\cr\cr
&=&{\Delta_0(\bfw_{k-1})\prod_{j=1}^{k-1}(w_k-w_j)\over w_k^{r-k}\cdot w_k^{k-1}(w_1\cdots w_{k-1})^{r-1}}
\ovsig_+(\bfw_{k})
 \Delta_\blamb(\sigma_-(w_k)\sigma_-(\bfw_{k-1})H_{r-k}) \wb^{r-k}_0\cr\cr\cr
 &=&{\Delta_0(\bfw_k)\over (w_1\cdots w_{k-1}\cdot w_k)^{r-1}}\cdot \ovsig_+(\bfw_k)
 \Delta_\blamb(\sigma_-(\bfw_{k})H_{r-k}) \wb^{r-k}_0\cr\cr\cr
 &=&{\Delta_0(\bfw_k)\over (w_1\cdots w_{k-1}\cdot w_k)^{r-1}}\Gamma^*(\bfw_k)\wb^r_\blamb,
\end{eqnarray*}
as claimed.\qed

\section{The main Theorem and its declinations}\label{sec:sec7}

In this section we shall be concerned with the several declinations of the main theorem 
describing the $B_r$ representation of $gl(\bw^kV)$.  
\claim{} Let $gl(\bw^kV):=\bw^kV\otimes\bw^kV^*$ be  the restricted Lie algebra of 
endomorphisms of $\bw^kV$, with respect to the natural commutator. With the same 
notation as in, a basis of $\bw^kV\otimes \bw^kV^*$ is $\big(\wb^k_\bmu\otimes\wbet^k_\bnu\big)_{\bmu,\bnu\in\Pcal_k}$, i.e.
$$
\bw^kV\otimes \bw^kV^*=\bigoplus_{\bmu,\bnu\in\Pcal_k}\QQ\cdot\wb^k_\bmu\otimes\wbet^k_\bnu,
$$
where $\wbet^k_\bnu(\wb^k_\bmu)=\delta_{\bmu,\bnu}$.

Then the $gl(\bw^kV)$-module structure of $B_r$ is defined through the following equality holding in $\bw^rV$:
\be
\big(\wb^k_\bmu\otimes\wbet^k_\bnu \star \Delta_\blamb(H_r)\big)\wb^r_0=\wb^k_\bmu\w (\wbet^k_\bnu\lrcorner \wb^r_\blamb).
\ee
This action is very easy to describe in the case $k=r$, but it becomes trickier when $r-k>0$. To describe it we shall consider the generating function
$$
\Ecal(\bfz_k,\bfw_k^{-1}):=\sum_{\bmu,\bnu}\wb^k_\bmu\otimes \wbet^k_\bnu\cdot  s_\bmu(\bfz_k)s_\bnu(\bfw_k^{-1}):B_r\sra B_r[[\bfz_k,\bfw_k]][\bfz_k^{-1},\bfw_k^{-1}].
$$
Our main result will consist in the explicit description of $\Ecal(\bfz_k,\bfw_k^{-1})\Delta_\blamb(H_r)$ in case $k\leq r$ (because otherwise one would obtain the trivial null action), where $\Ecal(\bfz_k,\bfw_k^{-1})\Delta_\blamb(H_r)$ is such that
\be
\Big(\Ecal(\bfz_k,\bfw_k^{-1})\Delta_\blamb(H_r)\Big)\wb^r_0=\sum_{\bmu\in\Pcal_k}s_\bmu(\bfz_k)\wb^k_\bmu\w (s_\bnu(\bfw_k ^{-1})\wbet^k_\bnu\lrcorner \wb^r_\blamb),
\ee
and where $s_\bmu(\bfz_k)$ and $s_\bnu(\bfw_k^{-1})$ denote the Schur symmetric polynomials labeled by the respective partitions with respect to the variables $\bfz_k$ and $\bfw_k^{-1}$ respectively.

%A basis of $\bw^kV^*$ is $\wbet^k_\blamb$ and $\wbet^k_\bmu(\wb^k_\blamb)=\delta_ 
%{\blamb,\mu}$. Therefore
%$$
%gl(\bw^kV)=\bigoplus_{\blamb\bmu\in\Pcal_k}\QQ\cdot \Ecal^k_{\blamb\bmu},
%$$
%where $\Ecal^k_{\blamb\bmu}:=\wb^k_\bmu\otimes\wbet^k_\bnu$
%and $B_r$ represents $gl(\bw^kV)$ in the obvious way:
%$$
%(\Ecal^k_{\blamb\bmu}\Delta_\bnu(H_r))\wb^r_0=\wb^k_\blamb\w (\wbet^k_\bmu\lrcorner \wb^r_\blamb)
%$$
%Our goal is to compute the generating series
%$$
%\Ecal(\bfz_k,\bfw^{-1}_k)\Delta_\blamb(H_r):=\sum_{\bmu,\bnu}\Ecal^k_{\bmu,\bnu}\Delta_\blamb(H_r)s_\bmu(\bfz_k)s_\bnu(\bfw_k^{-1})
%$$
%defined via the equality
%$$
%\sum_{\bmu,\bnu}(\Ecal^k_{\bmu,\bnu}\Delta_\bnu(H_r))\wb^r_0=\sum_{\bmu,\bnu}(\Ecal^k_{\bmu\bnu}\wb^r_\blamb)s_\bmu(\bfz_k)s_\bnu(\bfw_k^{-1})=\sum_{\bmu,\bnu\in\Pcal_k}\Big(\wb^k_\bmu\w (\wbet^k_\bnu\lrcorner \wb^r_\blamb)\Big)s_\bmu(\bfz_k)s_\bnu(\bfw_k^{-1})
%$$

\bclm{\bf Lemma.}\label{lem:genbl} {\em The generating function of the basis $(\wb^k_\bmu)_{\bmu\in\Pcal_k}$ of $\bw^kV$ is:
$$
\sum_{\bmu\in\Pcal_k}s_\bmu(\bfz_k)\wb^k_\bmu=\sigma_+(\bfz_k)\wb^k_0.
$$
}
\eclm
\proof By exploiting the definition of the $B_k$-module structure of $\bw^kV$, we have
\begin{center}
	\begin{tabular}{rclrr}
$\displaystyle{\sum_{\bmu\in\Pcal_k}}s_\bmu(\bfz_k)\wb^k_\bmu$&$=$&$\displaystyle{\sum_{\bmu\in\Pcal_k}}s_\bmu(\bfz_k)\Delta_\bmu(H_k)\wb^k_0$&&(using the $B_k$-module \cr
&&&&structure of $\bw^kV$ )\cr\cr
%&=&\sum_{\bmu\in\Pcal_k}s_\bmu(\bfz_k)\Delta_\bmu(H_k)\wb^k_0\cr\cr
&$=$&$\left(\displaystyle{\sum_{\bmu\in\Pcal_k}}s_\bmu(\bfz_k)\Delta_\bmu(H_k)\right)\wb^k_0$\cr\cr\cr
&$=$&$\displaystyle{\prod_{j=1}^k}\big(1+h_1z_j+h_2z_j^2+\cdots \big)\wb^k_0$&\cr\cr
&$=$&$\displaystyle{1\over E_k(z_1)}\cdot \displaystyle{1\over E_k(z_2)}\cdots{1\over E_k(z_k)}\wb^k_0$&& (By Cauchy formula as\cr
&&&& in \cite[Proposition 2, (iii)]{Fulyoung})\cr\cr
&$=$&$\sigma_+(\bfz_k)\wb^k_0$. &&\hskip100pt \qed
\end{tabular}
\end{center}
%
%\begin{eqnarray*}
%\sum_{\bmu\in\Pcal_k}s_\bmu\wb^k_\bmu&=&\sum_{\bmu\in\Pcal_k}s_\bmu(\bfz_k)\Delta_\bmu(H_k)\wb^k_0\cr\cr\cr
%&=&\sum_{\bmu\in\Pcal_k}s_\bmu(\bfz_k)\Delta_\bmu(H_k)\wb^k_0\cr\cr
%&=&\left(\sum_{\bmu\in\Pcal_k}s_\bmu(\bfz_k)\Delta_\bmu(H_k)\right)\wb^k_0\cr\cr\cr
%&=&\prod_{j=1}^k(1+h_1z_j+h_zz_j^2+h_3z_j^3+\cdots)\wb^k_0\cr\cr
%&=&{1\over E_k(z_1)}\cdot {1\over E_k(z_2)}\cdots{1\over E_k(z_k)}\wb^k_0\cr\cr
%&=&\sigma_+(z)\wb^k_0.\hskip 300pt\qed
%\end{eqnarray*}

\bclm{\bf Lemma.}\label{lem:genbet} {\em The generating function of the basis elements $\bw^kV^*$ is:	
\be
\sum_{\bnu\in\Pcal_k}s_\bnu(\bfw_k^{-1})\wbet^k_\bnu={\prod_{j=1}^kw_j^{k-1}\over \Delta_0(\bfw_k)}\cdot \bbeta(w_k^{-1})\w\cdots\w \bbeta(w_1^{-1}). \label{eq7:54}
\ee
}
\eclm
\proof
The one we propose consists in expanding the wedge product of the generating series of the basis $(\beta_j)_{j\geq 0}$ of $V^*$:
\begin{eqnarray}
&&\bbeta(w_k^{-1})\w\cdots\w \bbeta(w_1^{-1})\label{eq7:56}\\ \cr
&=&\sum_{\bnu\in\Pcal_k}\sum_{\tau\in S_k}\sgn(\tau)\beta_{k-\tau(1)+\nu_{\tau(1)}}\w\beta_{k-\tau(2)+\nu_{\tau(2)}} \cdots\w \beta_{k-\tau(k)+\nu_{\tau(k)}}\cdot\cr\cr\cr
&&\cdot w_k^{-k+\tau(1)-\nu_{\tau(1)}}w_{k-1}^{-k+\tau(2)-\nu_{\tau(2)}}\cdots w_1^{-k+\tau(k)-\nu_{\tau(k)}}=
\sum_{\bnu\in\Pcal_k}\wbet^k_\bnu\Delta_\bnu(\bfw_k^{-1})\cr\cr\cr
&=&\sum_{\bnu\in\Pcal_k}\wbet^k_\bnu \cdot s_\bnu(\bfw_k^{-1})\Delta_0(\bfw_k^{-1})=\sum_{\bnu\in\Pcal_k}s_\bnu(\bfw_k^{-1})\wbet^k_\bnu\cdot {\Delta_0(\bfw_k)\over w_k^{k-1}\cdots w_1^{k-1}},
\end{eqnarray}
whence the claim, obtained by multiplying both \eqref{eq7:54}  and \eqref{eq7:56} by $(\prod_{j=1}^kw_j^k)/\Delta_0(\bfw_k)$.\qed

\bclm{\bf Lemma.} {\em For all $u\in \bw^rV[[\bfw_k,\bfw_k^{-1}]$
\be
\sum_{\bmu}(\wb^k_\bmu\w u) s_\bmu(\bfz_k)=\ovsig_+(\bfz_k)\wb^k_0\w u={\prod_{j=1}^kz_j^r}\Gamma(\bfz_k)u.
\ee
}
\eclm
\proof We specified that the equality holds in $\bw^rV\llb\bfw_k,\bfw_k^{-1}\rrb$ to  emphasize the supposed $\QQ\llb \bfw_k,\bfw_k^{-1}\rrb$ linearity of the Schubert derivation. This said, it is not restrictive to assume that $u$ is a basis element $\wb^r_\blamb$ of $\bw^rV$.
The basic remark is that
\begin{eqnarray*}
{1\over E_k(z_1)}\cdot {1\over E_k(z_2)}\cdots{1\over E_k(z_k)}&=&\prod_{j=1}^k(1+h_1z_j+h_zz_j^2+h_3z_j^3+\cdots)\cr\cr
&=&\sum_{\bmu\in\Pcal_k}s_\bmu(\bfz_k)\Delta_\blamb(H_k),
\end{eqnarray*}
where in the last equality we used one of the declination of the celebrated Cauchy formula, as in \cite[Proposition 2, (iii)]{Fulyoung}, already used in the proof of Lemma \ref{lem:genbl}. Therefore:
\begin{eqnarray*}
\sum_\bmu\wb^k_\mu\ s_\bmu(\bfz_k)&=&\sum_\bmu\big(\Delta_\bmu(H_k) s_\bmu(\bfz_k)\big)\wb^k_0\cr\cr
&=&{1\over E_k(z_1)\cdots E_k(z_k)}\wb^k_0=\sigma_+(\bfz_k)\wb^k_0,
\end{eqnarray*}
where in the last equality we have repeatedly used the module structure of $\bw^rV$ over $B_r$.
Then:
$$
\sum_\bmu\wb^k_\mu\ s_\bmu(\bfz_k)\w u=\sigma_+(\bfz_k)\wb^k_0\w u,
$$
and the result now follows from Theorem~\ref{thm:thm61}.\qed

%We need one more preliminary result:

%\bclm{\bf Lemma.}
%
%\eclm
%\proof 
%Let us notice that
%
%Therefore
%\begin{eqnarray*}
%\sum_{\bnu\in\Pcal_k}(\wbet^k_\bnu \lrcorner \wb^r_\blamb)s_\bnu(\bfw_k^{-1}){\Delta_0(\bfw_k)\over w_k^{k-1}\cdots w_1^{k-1}}={\Delta_0(\bfw_k)\over w_1^{r-1}\cdots w_1^{r-1}}\Gamma^*(\bfw_k)\wb^r_\blamb
%\end{eqnarray*}	
%from which (\ref{eq7:54}), after clearing the common factors on both sides.\qed

We can finally express the action of the generating function $\Ecal_{\bmu,\bnu}(\bfz_k,\bfw_k^{-1})$ on a basis element of $B_r$.

\bclm{\bf Theorem (First Version).}\label{thm:mainthm} {\em For all $\blamb\in\Pcal_r${\em :}
\begin{eqnarray}
(\Ecal(\bfz_k,\bfw_k^{-1})\Delta_\blamb(H_r))\wb^r_0 &=&\prod_{j=1}^k\left({z_j\over w_j}\right)^{\hskip-4pt{r-k}}\Gamma(\bfz_k)\Gamma^*(\bfw_k)\wb^r_\blamb.
\end{eqnarray}	
}
\eclm
\proof
We have
\begin{center}
\begin{tabular}{rrll}
$\Ecal(\bfz_k,\bfw_k^{-1})\wb^r_\blamb$&$\hskip-7pt=$&$\hskip-7pt\displaystyle{\sum_{\bmu\in\Pcal_k}}s_\bmu(\bfz_k)\wb^k_\bmu\w \sum_{\bnu\in\Pcal_k}s_\bnu(\bfw_k)\Big(\wbet^k_\bnu\lrcorner \wb^r_\blamb\Big)$&\hskip-10pt (definition of \cr
&&&$\Ecal(\bfz_k,\bfw_k^{-1})$)\cr\cr
&$=$&$\sigma_+(\bfz_k)\wb^k_0\w \displaystyle{\prod_{j=1}^kw_j^{k-1}\over \Delta_0(\bfw_k)}\cdot \bbeta(w_k^{-1})\w\cdots\w \bbeta(w_1^{-1})\lrcorner \wb^r_\blamb$&(Lemmas\cr &&&\ref{lem:genbl} and \ref{lem:genbet})\cr\cr
&$=$&$\sigma_+(\bfz_k)\wb^k_0\w \displaystyle{1\over \prod_{j=1}^kw_j^{r-k}}\Gamma^*(\bfw_k)\wb^r_\blamb$.&\hskip-10pt (Theorem \ref{thm:thm71})
\end{tabular}
\end{center}
Now we use the fact that $\Gamma^*(\bfw_k)\wb^r_\blamb$ is a $\QQ[[\bfw_k,\bfw_k^{-1}]$-linear combination of $\wb^{r-k}_\bmu$ and then, by invoking Theorem~\ref{thm:thm61}, applied to last equality above,  becomes:
$$
\prod_{j=1}^{k}\left({z_j\over w_j}\right)^{r-k}\Gamma(\bfz_k)\Gamma^*(\bfw_k)\wb^r_\blamb,
$$
as announced.\qed

\bclm{\bf Corollary.} \label{cor:cor75} {\em If $r-k\geq \ell(\blamb)$ then
$$
\Gamma(\bfz_k)\Gamma^*(\bfw_k)\wb^r_\blamb=\sigma_+(\bfz_k)\ovsig_-(\bfz_k)\ovsig_+(\bfw_k)\sigma_-(\bfw_k)\wb^r_\blamb.
$$
}
\eclm
\proof
First of all notice that for  every $\blamb\in\Pcal_r$, it turns out that
$$
\ovsig_+(\bfw_k)\sigma_-(\bfw_k)\wb^{r-k}_\blamb=\sum_{\bmu\in\Pcal_r}a_\bmu(\bfw_k,\bfw_k^{-1})\wb^{r-k}_\bmu,
$$
where $a_\bmu(\bfw_k,\bfw_k^{-1})\in\QQ[\bfw_k,\bfw_k^{-1}]$. 
Then we have:
\begin{eqnarray*}
\Gamma(\bfz_k)\Gamma^*(\bfw_k)\wb^r_\blamb&=&\Gamma(\bfz_k)\ovsig_+(\bfw_k)\sigma_-(\bfw_k)\wb^{r-k}_\blamb\cr\cr
&=&\sum_{\bmu\in\Pcal_r}a_\bmu(\bfw_k,\bfw_k^{-1})\Gamma(\bfz_k)\wb^{r-k}_\bmu\cr\cr
&=&\sum_{\bmu\in\Pcal_r}a_\bmu(\bfw_k,\bfw_k^{-1})\wb^r_\blamb\cr\cr
&=&\sigma_+(\bfz_k)\ovsig_-(\bfw_k)\sum_{\bmu\in\Pcal_r}a_\bmu(\bfw_k,\bfw_k^{-1})\wb^r_\bmu\cr\cr
&=&\sigma_+(\bfz_k)\ovsig_-(\bfw_k)\sigma_+(\bfz_k)\ovsig_-(\bfw_k)\wb^r_\blamb.
\end{eqnarray*}
\claim{\bf Remark.}
If $\ell(\blamb)>r-k$, Corollary~\ref{cor:cor75} fails. We have however the following uniform way to compute $\Gamma^*(\bfw_k)\wb^r_\blamb$.

\bclm{\bf Proposition.}\label{prop:78} {\em For all $\blamb\in\Pcal_r$ and $k,r\geq 0${\rm :}
$$
\Gamma^*(\bfw_k)\wb^r_\blamb=\begin{vmatrix}\displaystyle{1\over w_1^{r-1+\lambda_1}}&\displaystyle{1\over w_1^{r-2+\lambda_2}}&\cdots&\displaystyle{1\over w_1^{\lambda_r}}
\cr
\vdots&\vdots&\ddots\cr\cr
\displaystyle{1\over w_k^{r-1+\lambda_1}}&\displaystyle{1\over w_k^{r-2+\lambda_2}}&\cdots&\displaystyle{1\over w_k^{\lambda_r}}
\cr\cr\cr
b_{r-1+\lambda_1}&b_{r-2+\lambda_2}&\cdots&b_{\lambda_r}
\end{vmatrix}\in \bw^{r-k}V,
$$
}
\eclm
using the same notation as in~\ref{sec4:sec41}.

\proof
By Theorem~\ref{thm:thm71}
\bigskip
\begin{eqnarray}
\Gamma^*(\bfw_k)\wb^r_\blamb&=&{(z_1\cdots z_k)^{r-1}\over \Delta_0(\bfw_k)}(\bbeta(w_k^{-1})\w \bbeta(w_{k-1}^{-1})\w\cdots\w \bbeta(w_1 ^{-1}))\lrcorner\wb^r_\blamb\cr\cr
&=&{(z_1\cdots z_k)^{r-1}\over \Delta_0(\bfw_k)}\begin{vmatrix}\displaystyle{1\over w_1^{r-1+\lambda_1}}&\displaystyle{1\over w_1^{r-2+\lambda_2}}&\cdots&\displaystyle{1\over w_1^{\lambda_r}}
\cr
\vdots&\vdots&\ddots\cr\cr
\displaystyle{1\over w_k^{r-1+\lambda_1}}&\displaystyle{1\over w_k^{r-2+\lambda_2}}&\cdots&\displaystyle{1\over w_k^{\lambda_r}}
\cr\cr\cr
b_{r-1+\lambda_1}&b_{r-2+\lambda_2}&\cdots&b_{\lambda_r}
\end{vmatrix}.\label{eq7:for63}
\end{eqnarray}
\qed

\noindent
Let us agree that
\smallskip
$$
w_k^{r-k}\cdots w_2^{r-2}w_1^{r-1}\Delta_\blamb(\bfw_k^{-1},H_{r-k})\wb^r_0=\bbeta(w_k^{-1})\w\cdots\w \bbeta(w_1^{-1}))\lrcorner\wb^r_\blamb,
$$

\noindent
\smallskip
defines $\Delta_\blamb(\bfw_k,H_{r-k})\in B_r[w^{-1}]$. The expansion of \eqref{eq7:for63} as a linear combinations of basis elements of $\bw^{r-k}V$, Giambelli's formula \eqref{eq:gmbsh} and the expansion rule of a determinant, easily imply that
\be
\Delta_\blamb(\bfw_k^{-1},H_{r-k})=\begin{vmatrix}\displaystyle{1\over w_1^{\lambda_1}}&\displaystyle{1\over w_1^{\lambda_2-1}}&\cdots&\displaystyle{1\over w_1^{\lambda_r-r+1}}
\cr
\vdots&\vdots&\ddots\cr\cr
\displaystyle{1\over w_k^{\lambda_1+k-1}}&\displaystyle{1\over w_k^{\lambda_2+k-2}}&\cdots&\displaystyle{1\over w_k^{\lambda_r+k-r}}
\cr\cr
h_{\lambda_1+k}&h_{\lambda_2+k+1}&\ldots&h_{\lambda_r+k+r-1}\cr
\vdots&\vdots&\ddots&\vdots\cr
h_{\lambda_1+r-1}&h_{\lambda_2+r-2}&\ldots&h_{\lambda_r}\label{eq:forfin}
\end{vmatrix}.
\ee

\noindent
This enables to state  a second version of \ref{thm:mainthm}, which works well for practical purposes and generalises \cite[Main Theorem 4.3]{gln}:
\bclm{\bf Theorem (second Version).}\label{thm:them79}
{\em The following equality holds {\em :}
\begin{eqnarray}
&&\Ecal(\bfz_k,\bfw_k^{-1})\Delta_\blamb(H_r)=\prod_{j=1}^k{\left(z_j\over w_j\right)^{r-k}}{1\over E_r(z_j)}\Delta_\blamb(\bfw_k^{-1},\ovsig_-(z)H_r)\cr 
&\hskip-20pt=\,\,\,&\hskip-15pt \prod_{j=1}^k{\left(z_j\over w_j\right)^{r-k}}{1\over E_r(z_j)}\begin{vmatrix}\displaystyle{1\over w_1^{\lambda_1}}&\displaystyle{1\over w_1^{\lambda_2-1}}&\cdots&\displaystyle{1\over w_1^{\lambda_r-r+1}}
\cr
\vdots&\vdots&\ddots\cr\cr
\displaystyle{1\over w_k^{\lambda_1+k-1}}&\displaystyle{1\over w_k^{\lambda_2+k-2}}&\cdots&\displaystyle{1\over w_k^{\lambda_r+k-r}}
\cr\cr
\ovsig_-(\bfz_k)h_{\lambda_1+k}&\ovsig_-(\bfz_k)h_{\lambda_2+k+1}&\ldots&\ovsig_-(\bfz_k)h_{\lambda_r+k+r-1}\cr
\vdots&\vdots&\ddots&\vdots\cr
\ovsig_-(\bfz_k)h_{\lambda_1+r-1}&\ovsig_-(\bfz_k)h_{\lambda_2+r-2}&\ldots&\ovsig_-(\bfz_k)h_{\lambda_r}
\end{vmatrix},
\end{eqnarray}
}
\eclm
where 
\be
\ovsig_-(\bfz_k)h_j=h_j-e_1(\bfz_k^{-1})h_{j-1}+\cdots+(-1)^ke_k(\bfz_k^{-1})h_{j-k}.\label{eq7:s-bzk}
\ee
\proof
By Theorem \ref{thm:mainthm} we have:

\begin{center}
\begin{tabular}{llcr}
$\Ecal(\bfz_k,\bfw_k^{-1})\Delta_\blamb(H_r)$&$=$&$\displaystyle{\prod_{j=1}^k}\left(\displaystyle{z_j\over w_j}\right)^{\hskip-4pt{r-k}}\Gamma(\bfz_k)\Gamma^*(\bfw_k)\Delta_\blamb(H_r)$&(Theorem \ref{thm:mainthm})\cr\cr
&$=$&$\displaystyle{\prod_{j=1}^k}\left(\displaystyle{z_j\over w_j}\right)^{\hskip-4pt{r-k}}\Gamma(\bfz_k)\Delta_\blamb(\bfw_k,H_{r-k})$&(Proposition~\ref{prop:78} \cr 
&&& and equation \eqref{eq:forfin})\cr\cr
&$=$&$\displaystyle{\prod_{j=1}^k\left({z_j\over w_j}\right)}^{\hskip-4pt{r-k}}\sigma_+(\bfz_k)\ovsig_-(\bfz_k)\Delta_\blamb(\bfw_k,H_{r})$&(Definition \cr
&&&\ref{def4:vertex}--\eqref{eq4:gamma} of $\Gamma(\bfz_k)$)\cr\cr
\cr\cr
&$=$&$\displaystyle{\prod_{j=1}^k\left({z_j\over w_j}\right)}^{\hskip-4pt{r-k}}{1\over E_r(z_j)}\cdot \Delta_\blamb(\bfw_k,\ovsig_-(\bfz_k)H_{r}) $&(Definition \ref{sec3:mosst}\cr &&&of the $B_r$-module\cr &&&structure of $\bw^rV$ and\cr
&&& Proposition~\ref{prop3:prop35})
\end{tabular}
\end{center}
as desired. Expression \eqref{eq7:s-bzk} for $\ovsig_-(\bfz_k)h_j$ is Lemma~\ref{lemma:53}.\qed
\claim{} Finally, let us define, as it is customary, new formal variable $(x_j)_{j\geq 1}$ through the equality:
$$
\exp\left(\sum_{j\geq 0}x_jz^j\right)={1\over E_r(z)}.
$$
In this case one can write
$$
\prod_{j=1}^k{1\over E_r(z_j)}=\exp\left(\sum_{j=0}x_jp_j(\bfz_k)\right),
$$
where $p_j(\bfz_k)=z_1^j+\cdots+z_k^j$ is the $i$-th power sum symmetric polynomial in $(z_1,\ldots,z_k)$ and where $x_i$ is precisely the $i$-th degree power sum in the $r$ universal roots $(y_1,\ldots,y_r)$ of the polynomial $E_r(z)$, i.e. $E_r(z)=\prod_{i=1}^r(1-y_jz)$ in the universal splitting $\QQ$-algebra for the polynomial $E_r(z)\in\QQ[z]$. This allows to shape our result in the form
\bclm{\bf Corollary.} {\em We have:
\be
\Ecal(\bfz_k,\bfw_k^{-1})\Delta_\blamb(H_r)=\prod_{j=1}^k{\left(z_j\over w_j\right)^{r-k}}\exp\left(\sum_{j=0}x_jp_j(\bfz_k)\right)\Delta_\blamb(\bfw_k,\ovsig_-(\bfz_k)H_{r-k}).\label{eq:oulast}
\ee
}
\eclm
\qed

Formula \eqref{eq:oulast} is easy to use for  practical computations of the $gl(\bw^kV)$--representation of $B_r$,  for those special case of $k$ and $r$ that everybody may possibly need.

\bigskip

\noindent{\bf Acknowledgments.} This work is an expansion of part of the Ph.D. thesis of the first author performed  during her hosting at the Department of Mathematical Sciences of Politecnico of Torino under the sponsorship of Ministry of Science of the Islamic Republic of Iran. The second and fourth authors stay enjoyed of  the full sponsorship of Finanziamento Diffuso della Ricerca, no. 53$\_$RBA17GATLET and Progetto di Eccellenza del Dipartimento di Scienze Matematiche, 2018--2022, no.
E11G18000350001. The project also benefitted the partial 
support of INDAM-GNSAGA, PRIN ``Geometria delle Variet\`a Algebriche''.

For discussions and criticisms we want to 
thank  primarily Parham Salehyan, Abbas Nasrollah Nejad, Inna Scherbak, Joachim Kock,  Piotr Pragacz and Andrea T. Ricolfi for many kind of assistance.

%\begin{eqnarray*}		
%		&=&\exp(\sum_{i\geq 1}x_ip_i(\bfz_k))\exp\left(\sum_{i\geq 1}{1\over i} \delta(\sigma_{-1}^i)p_i(\bfz_k)\right)\ovsig_+(\bfw_k)\Delta_\blamb(\bfw_k^{-1},\prod_{i=1}^k\ovsig_-(z_i)H_{r})\\
%			&=&\exp(\sum_{i\geq 1}x_ip_i(\bfz_k))\exp\left(\sum_{i\geq 1}{1\over i} \delta(\sigma_{-1}^i)p_i(\bfz_k)\right)\exp\left(-\sum_{i\geq 1}x_ip_i(\bfw_k)\right)\Delta_\blamb(\bfw_k^{-1},\prod_{i=1}^k\ovsig_-(z_i)H_{r})
%\end{eqnarray*}
%In other words $\wb^k_\bmu\otimes \wbet^k_\bnu$ is the coefficent of $\bfz_k^\bmu\bfw_k^{-\bnu}$ the expansion of
%$$
%{\Delta_0(\bfz_k)\cdot \Delta_0(\bfw_k^{-1})}\Ecal(\bfz_k,\bfw_k^{-1})\Delta_\blamb(\bfw_k^{-1},H_r).
%$$

%\newpage
\bibliographystyle{amsplain}
%\bibliography{ourrefs}
\providecommand{\bysame}{\leavevmode\hbox to3em{\hrulefill}\thinspace}
\providecommand{\MR}{\relax\ifhmode\unskip\space\fi MR }
\providecommand{\MRhref}[2]{%
  \href{http://www.ams.org/mathscinet-getitem?mr=#1}{#2}
}
\providecommand{\href}[2]{#2}

\medskip
\medskip

\parbox[t]{3in}{{\rm Ommolbanin~Behzad}\\
{\tt \href{mailto:behzad@iasbs.ac.ir}{behzad@iasbs.ac.ir}}\\
{\it Institute for Advanced Studies in Basic Sciences, Zanjan}\\
{\it IRAN}} \hspace{1.5cm}
\parbox[t]{3in}{{\rm Andr\'e Contiero}\\
{\tt \href{mailto:contiero@ufmg.br}{contiero@ufmg.br}}\\
{\it Universidade Federal de Minas Gerais}\\
{\it Belo Horizonte, MG}\\
{\it BRAZIL}}

\bigskip

\vspace{6 pt}

\parbox[t]{3in}{{\rm Letterio~Gatto}\\
{\tt \href{mailto:letterio.gatto@polito.it}{letterio.gatto@polito.it}}\\
{\it Dipartimento~di~Scienze~Matematiche}\\
{\it Politecnico di Torino}\\
{\it ITALY}} \hspace{1.5cm} 
\parbox[t]{3in}{{\rm Renato Vidal Martins}\\
{\tt \href{mailto:vidalmartins@ufmg.br}{vidalmartins@ufmg.br}}\\
{\it Universidade Federal de Minas Gerais}\\
{\it Belo Horizonte, MG}\\
{\it BRAZIL}}

\end{document}